%% file: article_color.tex
\documentclass[onefignum,onetabnum]{siamonline190516}

\input{shared}

\ifpdf
    \hypersetup{
        pdftitle={Learning Optical Flow for Fast MRI Reconstruction},
        pdfauthor={T. Schmoderer, A.I Aviles-Rivero, V. Corona, N. Debroux and C-B. Schönlieb}
    }
\fi

\usepackage{subcaption}
\usepackage{graphicx}
\usepackage{comment}
\usepackage{booktabs}
\usepackage{multirow}
\usepackage{wrapfig}

\begin{document}

\maketitle

\begin{abstract}
Reconstructing high-quality magnetic resonance images (MRI) from undersampled raw data is of great interest from both technical and clinical point of views. To this date, however, it is still a mathematically and computationally challenging problem due to its severe ill-posedness, resulting from the highly undersampled data leading to significantly missing information. Whilst a number of techniques have been presented to improve image reconstruction, they only account for spatio-temporal regularisation, which shows its limitations in several relevant scenarios including dynamic data. In this work, we propose a new mathematical model for the reconstruction of high-quality medical MRI from few measurements. Our proposed approach combines - \textit{in a multi-task and hybrid model} - the traditional compressed sensing formulation for the reconstruction of dynamic MRI with motion compensation by learning an optical flow approximation. More precisely, we propose to encode the dynamics in the form of an optical flow model that is sparsely represented over a learned dictionary. This has the  advantage that ground truth data is not required in the training of the optical flow term.  Furthermore, we present an efficient optimisation scheme to tackle the non-convex problem based on an alternating splitting method. We demonstrate the potentials of our approach through an extensive set of  numerical results using different datasets and acceleration factors. Our combined approach reaches and outperforms several state of the art techniques. Finally, we show the ability of our technique to transfer phantom based knowledge to real datasets.
\end{abstract}

\begin{keywords}
MRI Reconstruction, Optical Flow,  Dictionary Learning, Compressed Sensing, Multi-Task Model
\end{keywords}

\begin{AMS}
94A08, 65K05, 68Q32
\end{AMS}

\section{Introduction}

\input{sections/0_introduction.tex}
\section{MRIR-DLMC for Fast MRI: Proposed Model}\label{sec:model}
\input{sections/2_model.tex}

\section{MRIR-DLMC for Fast MRI: Numerical Realisation}\label{sec:alg}
\input{sections/3_algorithm.tex}

\section{Numerical Experiments}\label{sec:experiments}
\input{sections/4_numericals.tex}


\section{Conclusion}\label{sec:conclusions}
\input{sections/5_conclusion}

\appendix
\input{sections/6_appendices}

\section*{Acknowledgments}AIAR  gratefully acknowledges the financial support of the CMIH and CCIMI University of Cambridge. VC acknowledges the financial support of the Cancer Research UK. CBS acknowledges support from the Leverhulme Trust project on 'Breaking the non-convexity barrier', EPSRC grant Nr. EP/M00483X/1, EP/S026045/1, the EPSRC Centre Nr. EP/N014588/1,  Wellcome Trust project "All-in-one cancer imaging", the RISE projects CHiPS and NoMADS, the Cantab Capital Institute for the Mathematics of Information and the Alan Turing Institute.

\bibliographystyle{siamplain}
\bibliography{references}

\end{document}

%% file: shared.tex
\usepackage{lipsum}
\usepackage{amsfonts}
\usepackage{graphicx}
\usepackage{epstopdf}
\usepackage{algorithmic}
\ifpdf
  \DeclareGraphicsExtensions{.eps,.pdf,.png,.jpg}
\else
  \DeclareGraphicsExtensions{.eps}
\fi

\usepackage{enumitem}
\setlist[enumerate]{leftmargin=.5in}
\setlist[itemize]{leftmargin=.5in}

\newsiamremark{remark}{Remark}
\newsiamremark{hypothesis}{Hypothesis}
\crefname{hypothesis}{Hypothesis}{Hypotheses}
\newsiamthm{claim}{Claim}

\headers{MRIR-DLMC}{T. Schmoderer, A.I Aviles-Rivero, V. Corona, N. Debroux and C.-B. Sch\"onlieb}

\title{Learning Optical Flow for Fast MRI Reconstruction\thanks{Submitted to the editors on April 22, 2020.
}}

\author{
Timoth\'ee Schmoderer\thanks{Normandie Univ, Institut National des Sciences Appliqu\'ees de Rouen, Laboratory of Mathematics, 76000 Rouen, France 
  (\email{timothee.schmoderer@insa-rouen.fr}).}
\and Angelica I Aviles-Rivero\thanks{Department of Pure Mathematics and Mathematical Statistics, University of Cambridge, UK (\email{ai323@cam.ac.uk}).}
\and Veronica Corona\thanks{Department of Applied Mathematics and Theoretical Physics (DAMTP), Centre for Mathematical Sciences, University of Cambridge, Wilberforce Road, Cambridge CB3 OWA, UK
  (\email{vc324@cam.ac.uk},\email{cbs31@cam.ac.uk}).}
\and No\'emie Debroux\thanks{Universit\'e Clermont Auvergne, CNRS, SIGMA Clermont, Institut Pascal, F-63000 Clermont-Ferrand, France
  (\email{noemie.debroux@uca.fr}).}
  \and Carola-Bibiane Sch\"onlieb\footnotemark[4]}

\usepackage{amsopn}


%% file: sections/0_introduction.tex
The problem of reconstructing high quality images, from dynamic Magnetic Resonance (MR) measurements, whilst reducing the inherent involuntary motion during the acquisition is a central topic in the community. Although it has been widely explored, it is still considered an open problem that we address in this work. This central topic in MRI comes as a consequence of the long acquisition time, which makes the image formation highly sensitive to motion leading to image degradation and compromising the expert interpretation \cite{sachs1995diminishing,Zaitsev::2015}. 

There are several attempts in the body of literature, to improve the reconstruction under inherent motion. A set of solutions to reduce motion artefacts includes breath-holding techniques \cite{plathow2006assessment,Zaitsev::2015,gdaniec2014robust} and gating strategies \cite{Setser2000,Plein2001,Kaji2001,GEORGE2006924,JIANG2006141,Zaitsev::2015}.
The former rely on the  patients' ability to hold their breath. However, the time needed to form an image is much longer than the average time a human being is capable of holding its breath
~\cite{plathow2006assessment}. The latter aim at tracking either the breathing or the cardiac cycles via external sensors. However, the co-registration of these signals to the image ones is still challenging ~\cite{HOISAK2006339,Zaitsev::2015} and these techniques are effective only for perpetual motions disregarding any other involuntary motion or arrhythmia \cite{Zaitsev::2015} and thus are only partially accurate.

\textit{Fast MRI reconstruction.} The idea of accelerating the MRI acquisition by using under-sampling strategies comes as a complement 
to the aforementioned techniques. This option has 
gained great interest since the seminal paper \cite{lustig2007} in which Compressed Sensing techniques \cite{donoho2006compressed, candes2006stable} were successfully adapted to the MRI reconstruction problem. 
CS relies on the property of MR images to be sparse in a transformed domain. Examples are  \cite{LustigktSPARSE::2006,Otazo::2010,LustigSPIRiT::2010,lingala::2013,Liang::2007,gamper2008compressed,jung2009k,Lingala::2011}, in which authors applied several known or learnt transformations alone or in combination with parallel imaging. Extending the idea of CS, another set of works has been devoted to 
low-rank matrix completion by taking advantage of the highly redundant information present in dynamic MRI sequences. For example, the work of Liang~\cite{Liang::2007}, in which a singular value decomposition (SVD) is used to compute the temporal basis functions, has paved the way to techniques relying on Principal Component Analysis (PCA) to reconstruct from under-sampled measurements \cite{Pedersen::2009,Feng::2013,Velikina::2015,arif2019accelerated}. Other works have proposed to incorporate local low rank constraints in small regions to reduce the computational load \cite{Trzasko::2011,Zhang::2015,Miao::2016}. However they come at the expense of block artefacts \cite{Saucedo::2017}.

A body of research has proposed to combine sparse and low-rank constraints e.g.~\cite{Lingala::2011,Majumdar::2012,Zhao::2012}. 
Decomposing the reconstruction into a sparse component and a low rank one is known as both Robust Principal Component Analysis (RPCA) \cite{Candes::2011} and L+S decomposition \cite{Gao::2012,Tremoulheac::2014,Otazo::2013,Otazo::2015}. Even though the efficiency of these techniques has been proved in the context of dynamic MRI reconstruction, the incoherence, in which one needs to choose an irregular sampling pattern such that it does not replicate image features, between the low rank and sparse representations for robust separation of background and dynamic components required for the L+S decomposition is not always fulfilled. The success of these models is promising but there is still room for improvement as they do not take into account the explicit motion exhibited during the dynamic acquisition. Thus reconstruction techniques also considering this motion have emerged to improve the quality of the reconstructed sequence in a dynamic setting.

\textit{Motion Estimation.} There are two popular strategies for motion estimation based either on registration \cite{sotiras2013} or optical flow \cite{fortun2015}. In this work, we focus on the latter. Optical flow relies on the  brightness constancy assumption, and consists in estimating the velocities of movement of brightness patterns in an image. Optical flow is a very ill-posed problem due to the aperture problem.  
First, variational approaches were introduced to solve the OF problem, starting with the seminal work of Horn and Schunck \cite{horn1981}.
In that work, authors introduced an $L^2$ fidelity term coming from the linearised brightness constancy assumption and an $L^2$ regularity term on the flow fields to ensure smoothness and make the problem well-posed. Since then, there has been a large body of works trying to improve the quality and the accuracy of the obtained flows in a variational setting \cite{fortun2015,perez2013} (and references therein).

More recently, several authors have explored the optical flow problem using tools drawn from machine learning. Astonishing results have been reported using deep learning either for normal scene conditions e.g.\cite{sun2018pwc,dosovitskiy2015flownet,ilg2017flownet} or extreme scenarios e.g.\cite{li2019rainflow}. We remark that our purpose in this work is different from those works. Whilst they aim to improve the optical flow estimation itself, \textit{our goal is instead to use OF as a proxy task to improve the MRI reconstruction.} In particular, we are interested in learning the optical flow using sparse coding. Our motivation to follow this philosophy, is given by the promising results reported in image restoration e.g~\cite{elad2006image,mairal2009non}. The central idea of sparse coding is to find a dictionary such that each measurement can be well-approximated by a sparse linear combination of atoms (i.e. elements of a dictionary).

There have been several attempts to learn a sparse optical flow. The works of that~\cite{timofte2015sparse,shen2010sparsity}  used known basis functions whilst authors in~\cite{jia2011} learnt an overcomplete dictionary coming either from ground truth set or from classic variational estimation \cite{gibson2016sparsity,wulff2015efficient}. The major drawback of these techniques is the lack of convergence guarantee but this has been overcome by the work of Bao~\cite{bao2015dictionary}. In that work, authors introduced a globally convergent algorithm to solve generic dictionary learning and sparse coding problems. Whilst sparse coding has been widely explored for natural images, sparse optical flow models have yet to be tackled in the medical domain and only few works have been reported e.g.~\cite{ouzir2018robust}. 

Although previous approaches have shown potential results on both MRI reconstruction and Optical Flow estimation separately, more recently, a set of algorithmic approaches known as multi-tasking or joint models have demonstrated that the reconstruction can be further improved when a related task is included in the reconstruction loss. We remark that our ultimate goal is to improve the image reconstruction for dynamic MRI sequences. To this purpose, we include motion estimation as a proxy task for improving MRI reconstruction.

\textit{Joint Models.} Following this perspective, different works have been reported using so-called joint models, in which  the main idea is to intertwine several related tasks into a single unified functional. The potential in terms of image quality improvement whilst reducing error propagation has been demonstrated in~\cite{aviles2018, liu2019rethinking,caruana1997multitask,Zhao::2019}.
There are several works that follow this idea either for reconstructing a single image e.g.\cite{Odille::2008,Corona::2019,weller2019motion,Prieto::2007,haskell2018targeted} or a full image sequence e.g.~\cite{Lingala::2015,Odille::2016,Odille::2010,burger2018,Zhang::2018,Zhao::2019}. 
Some of them rely on prior and external knowledge to generate motion models from sensors, navigators, etc. \cite{Odille::2016} or driven by previously computed eigenmodes of cardiac and respiratory motions \cite{Odille::2010}. Iterative schemes, with two coupled inverse problems but not derived from the same original optimisation problem, 
also aim at computing simultaneously the motion map and the reconstructed images e.g.~\cite{royuela2016nonrigid,Royuela::2017,filipovic2011motion,Asif::2013,Rank::2017}. These models alternate between the resolution of an independent motion estimation problem, and the resolution of a motion informed reconstruction of MR sequences.

Authors of that \cite{Lingala::2015} proposed the self-contained Deformation Corrected-Compressed Sensing (DC-CS) model relying on CS principles along with a Demons algorithm to estimate the motion in a unified framework. More recently, Zhang et al. in~\cite{Zhang::2018} extended this model to the concept of Blind Compressed Sensing, in which the transformation leading to sparsity is also estimated. In particular, our proposed approach follows the same perspective as the works of \cite{aviles2018,Zhao::2019,burger2018}, in which CS and OF principles are combined. We emphasise that the main difference stands in the optical flow formulation. That is- we propose a new optimisation model, in which the flow estimation is encoded in a sparse representation over a learnt dictionary. 

In this work, we propose a new multi-task framework for reconstructing high quality MR images from reduced measurements. Our approach seeks to improve the reconstruction quality through a better optical flow estimation using a learning-based approach. More precisely, we plug-in to the reconstruction energy a learnt sparse optical flow loss. Our optical flow model is computed as a sparse linear combination of basis functions from a learnt dictionary. Whilst these are important part of our work, our contributions are as follows:

\begin{itemize}
    \item We introduce a unified joint and hybrid (i.e. a model that combines variational and learning strategies) model for fast MRI reconstruction, which we called \textbf{MRI R}econstruction with  \textbf{D}ictionary  \textbf{L}earnt  \textbf{M}otion  \textbf{C}ompensation (MRIR-DLMC), in which we highlight:
        \begin{itemize}
    	\item We propose a mathematical model that simultaneously performs the MRI reconstruction using Compressed Sensing principles and motion estimation via a TV-$L^1$-based dictionary learning approach for Optical Flow estimation.
    	\item We propose a tractable and efficient numerical scheme to solve the initial optimisation  problem.
    	\item We show that our proposed technique is able to do transfer knowledge from phantom to realistic data. 
        \end{itemize}
    \item We extensively evaluate our approach using several datasets and acceleration factors. From the results,
     \begin{itemize}
        \item We show that integrating motion information in the Compressed Sensing MRI reconstruction model improves the overall quality of the reconstructed sequence in comparison to state of the art MRI reconstruction techniques.
        \item  We demonstrate that the quality of the motion model is linked to the quality of the reconstruction, and that learned parametrisation for the motion as the one proposed in this paper outperforms simple motion models with handcrafted regularisation on the velocity field. 
        \item We demonstrate generalisation properties by showing improved performance over four datasets and using several undersampling factors. 
      \end{itemize}
\end{itemize}

%% file: sections/2_model.tex
This section is devoted to the depiction of our joint and hybrid model for simultaneous MRI reconstruction and motion estimation called  \textbf{MRI} \textbf{R}econstruction with \textbf{D}ictionary \textbf{L}earnt \textbf{M}otion \textbf{C}ompensation \textbf{MRIR-DLMC}. The section is organised in three key parts: i) the CS MRI reconstruction scheme ii) our proposed approach for learning optical flow; and iii) our combined MRI reconstruction and motion estimation  model.

\subsection{Compressed Sensing MRI Reconstruction} \label{sub-sec:CS_reconstruction}
We consider a dynamic MRI setting, where we denote by $\mathbf{f}$ the time sequence of raw data collected in a spatial-frequency space ($\mathbf{k}$,~$t$-space). The associated forward model -- relating the MR measurements $\mathbf{f}$ with the meaningful visual data $\mathbf{m}$ in the image domain $\Omega$, being a connected bounded open subset of $\mathbb{R}^2$ -- reads:
\begin{align}
    \mathbf{f}(\mathbf{k},t)=\int_\Omega \mathbf{m}(\mathbf{x},t)\exp({-j\mathbf{k}^T\cdot\mathbf{x}})\,d\mathbf{x}+\eta(\mathbf{k},t),
\end{align}
\noindent
where $\mathbf{x}$ is the spatial coordinate and $\mathbf{k}$ the frequency variable. Moreover, $t$ is the temporal coordinate varying from $0$ (initial acquisition time) to $T$ (final acquisition time),  $\eta$ is the inherent acquisition noise which can be modelled as Gaussian noise for the MRI application and $\mathbf{m}$ is  the image sequence of a moving part of the body.

For fast MRI reconstruction, one seeks to reconstruct the images from a highly reduced  number of acquired measurements, which translates into $\mathbf{f}(\mathbf{k},t)$ being available only for a small amount of $\mathbf{k}\in \mathbb{R}^2$. To do this, we denote by $K$ the undersampled Fourier operator, that is to say the composition of a sampling mask also encoding coil sensitivities with the Fourier transform. Nevertheless, the resulting inverse reconstruction problem becomes highly ill-posed and can only be solved by introducing  a prior knowledge on the nature of the reconstructed images.


The MRI reconstruction from highly undersampled data, using compressed sensing principles, was first investigated in the seminal paper of Lustig~\cite{lustig2007}. 
The main idea of CS in MRI is to take advantage of the natural sparsity of images in some transform domains, and especially the sparsity of MRI in the Wavelet domain.  
One can use the $L^1$ norm as a sparsity measure for which the unconstrained CS reconstruction minimisation problem  reads: 

\begin{align}\label{eq:classical_mri_reconstruction}
    \inf_{\mathbf{m}}  E(\mathbf{m}) = \int_0^T \frac{1}{2}\left\|K\mathbf{m}- \mathbf{f}\right\|^2_{L^2(\mathbb{R}^2)} + \lambda_1 TV(\mathbf{m}) + \lambda_2\left\|\Psi \mathbf{m}\right\|_{L^1(\mathbb{R}^4)}\ dt,
\end{align}
\noindent
where $TV(\mathbf{m})$ denotes the spatial total variation to enforce sparse edges and $\Psi$ indicates the Wavelet transform. Moreover, $\lambda_1$ and $\lambda_2$ are non-negative tuning parameters. The first term is a matching criterion to enforce closeness between the reconstructed images $\mathbf{m}$ and the MR measurements $\mathbf{f}$, while the remaining terms regularise the images by adding prior knowledge and making the problem well-posed. In this work, we focus on the total variation combined with the Wavelet transform favouring better medical image reconstruction with sharp edges, but this can be replaced by any other transformation easily such as the Fourier transform or the discrete cosine transform.


Although the body of literature has shown the potentials of using the model described in~\eqref{eq:classical_mri_reconstruction}, 
it is not time dependent and therefore is not performing well enough in dynamic settings, where motion and temporal redundancy are involved. With the purpose to improve the reconstruction quality, 
one can account for the inherent motion in the scene via motion estimation that can be computed using Optical Flow (OF). In what follows, we describe our new approach to improve OF approximation via a sparse representation in a learnt dictionary, and show its connection with MRI reconstruction.




\subsection{Motion Estimation: Learning Optical Flow} \label{sub-sec:motion_estimation}
Optical flow (OF) consists in estimating the vector field of apparent velocities of brightness patterns between consecutive frames. Assuming small displacements, the essential brightness constancy assumption can be linearised and reads:

\begin{align}
    \label{eq:brightness_constancy}
    \frac{\partial \mathbf{m}}{\partial t} + \mathbf{u}\cdot\nabla \mathbf{m} = 0,
\end{align}

\noindent
where $\mathbf{m} = \mathbf{m}(\mathbf{x}, t)$ denotes the dynamic sequence of images with $\frac{\partial \mathbf{m}}{\partial t}$ being the temporal derivative of the image sequence, $\nabla \mathbf{m}$ the spatial gradient, and $\mathbf{u}=[u_x,u_y]^T$ the unknown motion field. To deal with the aperture problem, one can embed \cref{eq:brightness_constancy} in a  variational formulation which reads, in its generic form: 

\begin{align}
    \inf_{\mathbf{u}} E(\mathbf{u})= \int_0^T E_{d}(\mathbf{u})+\lambda E_{r}(\mathbf{u})\, dt,
\end{align}

\noindent
where $E_{d}$ is a fidelity term of the form $E_{d}(\mathbf{u}) = \phi( \frac{\partial \mathbf{m}}{\partial t} + \mathbf{u}\cdot\nabla \mathbf{m})$, and $E_{r}$ is a regularisation of the flow field to make the problem well-posed. In their seminal paper \cite{horn1981}, Horn and Schunck set $\phi(\cdot)=\|\cdot\|_{L^2(\Omega)}^2$ and $E_{r}(\cdot)=\|\nabla \cdot\|_{L^2(\Omega)}^2$. Moreover, in the work of ~\cite{perez2013} authors proposed  $\phi(\cdot)=\|\cdot\|_{L^1(\Omega)}$ and $E_{r}(\cdot)=TV(\cdot)$, which is known as the $TV-L^1$ formulation. They show that this formulation leads to  better and more accurate results than the ones in~\cite{horn1981}. In the remainder of the paper, we build upon the $TV-L^1$ model.

\begin{remark}
    The brightness constancy assumption usually is not perfectly fulfilled for real-world problems, e.g.~\cite{sun2008learning,werlberger2010motion}, and in particular, with medical data due to its inherent properties. However, we show that this assumption still enables us to get a good approximation of the observed motion, for our application at hand, to improve the MRI reconstructions. A relaxation of this assumption to further correct the motion estimation leading to even better reconstruction quality will be the object of future work. 
\end{remark}

Whilst the variational formulations for the OF estimation is limited by its own design, one can further improve the approximation of the motion scene by expressing the flow field as a sparse linear combination of basis functions in either off-the-shelf dictionaries \cite{timofte2015sparse,shen2010sparsity} or learnt ones \cite{jia2011,wulff2015efficient,gibson2016sparsity}. Inspired by these works, we seek to increase the accuracy of the OF estimation by encoding it in a sparse representation over a learnt over-complete dictionary. From now on, to comply with the dictionary learning and sparse representation philosophies, we consider a discrete spatial setting. The overcomplete flow dictionary $D$ is thus learnt from a partition \textcolor{black}{$\mathcal{P}$} of an initial flow estimation decomposed into small overlapping patches.  In here, we define the patch extracting operator as $R_{\mathbf{p}}$ at pixel $\mathbf{p}$. We follow the strategy from~\cite{jia2011}, in which the flow dictionary $D$ is simplified by decoupling the horizontal and vertical motions as $D=\begin{bmatrix} D_x  & 0\\  0 & D_y \end{bmatrix}$, where both $D_x$ and $D_y$ are composed of $N_d$ elements such that $N_d\times |\mathcal{P}|$ is larger than the image dimension. Formally, this sparse representation over a learnt dictionary is incorporated in the variational setting through this additional term:

\begin{align} \label{eq:dicLear1}
    E_{sparse}(\mathbf{u},D,\mathbf{a}) &= \sum_{\mathbf{p}\in\mathcal{P}}\|R_{\mathbf{p}}\mathbf{u}-D\mathbf{a}_{\mathbf{p}}\|_F^2+\tau \|\mathbf{a}_{\mathbf{p}}\|_1\\
\nonumber    &s.t.\, \|D_{x,i}\|_2\leq1,\, \|D_{y,i}\|_2\leq1,\, 1\leq i \leq N_d,
\end{align}

\noindent
where $\|\cdot\|_F$ is the Frobenius norm and $\|\cdot\|_1$ is the $l^1$ norm. The new variable $\mathbf{a}_{\mathbf{p}}=((a_x,a_y)_{\mathbf{p}})^T$ represents the sparse coefficients of the patch $R_{\mathbf{p}}u$ in $D$. In the model \eqref{eq:dicLear1}, the first term ensures the flow field is given by a linear combination of basis functions in the dictionary whilst the second one enforces sparsity of the coefficients. At the learning level, the model \eqref{eq:dicLear1} can be solved using either ground truth when available \cite{jia2011}, or from an initial approximation of the flows using a variational model \cite{wulff2015efficient,gibson2016sparsity}. The assumption of having ground truth is not realistic in the medical domain, and therefore we follow the second strategy. Motivated by the results in \cite{wulff2015efficient,gibson2016sparsity} showing better results with a dictionary learnt from a classic variational formulation than the initial ones, we propose to build our dictionary upon the $TV-L^1$ model, and thus consider the following discrete optical flow formulation:

\begin{align}
    \label{eq:optical_flow_energy}
    E_{OF}(\mathbf{u},D,\mathbf{a}) &= \int_0^T \lambda_3\left\|\frac{\partial \mathbf{m}}{\partial t} + \nabla \mathbf{m}\cdot \mathbf{u} \right\|_{1} + \lambda_4 TV(\mathbf{u})+\sum_{\mathbf{p}\in\mathcal{P}}\lambda_5\|R_{\mathbf{p}}\mathbf{u}-D\mathbf{a}_{\mathbf{p}}\|_F^2+\lambda_6 \|\mathbf{a}_{\mathbf{p}}\|_1\\
\nonumber    &s.t.\, \|D_{x,i}\|_2\leq1,\, \|D_{y,i}\|_2\leq1,\, 1\leq i \leq n_d,
\end{align}

\noindent
where $\lambda_3,\ \lambda_4,\ \lambda_5\text{ and } \lambda_6 >0$ are weighting parameters to balance the influence of each term. The constraints on $D$ ensure uniqueness of the solution. Given a reference optical flow $\mathbf{u}^{ref}$, the dictionary is learnt as the result of the following minimisation problem under constraints:

\begin{align}
    \label{eq:learn_dico}
    E(D,\mathbf{a}) &= \int_0^T \|R_{\mathbf{p}}\mathbf{u}^{ref}_{\mathbf{p}}-D\mathbf{a}_{\mathbf{p}}\|_F^2\ dt, \\ 
    &\quad \|D_{x,j}\|_2 \leq 1,\ \|D_{y,j}\|_2 \leq 1, \quad 1\leq j\leq N_d, \nonumber \\
    &\quad \|\mathbf{a}_{\mathbf{p}}\|_0\leq k_0 \nonumber \ \forall \mathbf{p}\in\mathcal{P},
\end{align}

\noindent
where the pseudo norm $l^0$ counts the non-zero elements of $\mathbf{a_p}$. The model \cref{eq:learn_dico} uses a constraint on $\mathbf{a}$ rather than the $l^1$ relaxation of the $l^0$ pseudo-norm as it shows better numerical performances. 


\subsection{Joint Model: Learning Optical Flow for MRI Reconstruction}\label{sub-sec:hybrid_model} 

The potentials of CS MRI reconstruction and optical flow for motion estimation have been demonstrated separately in the previous subsections. We now demonstrate the strong correlation between both tasks -- as motion estimation depends on the reconstruction quality, and the reconstruction accuracy can be substantially improved by a good approximation of the motion scene -- and show their mutual benefit by intertwining them.
This statement has been shown in recent works such as \cite{aviles2018,Aviles-Rivero::2018,burger2018} that have inspired our work, in which one seeks to solve a single optimisation model unifying both tasks in a pure variational setting.
More precisely, the two formulations for CS MRI reconstruction and TV-$L^1$ Optical Flow estimation can be cast in a unified optimisation problem as: 

\begin{align}\label{initial_variational_problem}
     \inf_{\mathbf{m},\mathbf{u}}&\int_0^T \frac{1}{2}\|K\mathbf{m} - \mathbf{f}\|_{L^2(\mathbb{R}^2)}^2 + \lambda_1 TV(\mathbf{m})^p + \lambda_2\|\psi \mathbf{m}\|_{L^1(\mathbb{R}^4)}^p+\lambda_3\left\|\frac{\partial \mathbf{m}}{\partial t}+\mathbf{u}\cdot\nabla\mathbf{m}\right\|_{L^1(\Omega)} \\
     \nonumber &+ \lambda_4TV(\mathbf{u})^q\,dt,   
\end{align}

For now, we have just added the sparsity constraint on the Wavelet domain of the reconstruction sequence in a purely variational setting compared to \cite{burger2018,aviles2018,Aviles-Rivero::2018}. The existence theorem given in \cite{burger2018} still holds.

\begin{theorem}[Existence of minimisers]
Let the assumptions of \cite[Theorem 3.1]{burger2018} be satisfied (the MRI reconstruction operator fulfilling the required properties). Then there exists a minimiser of \cref{initial_variational_problem} in the space 

\begin{align*}
    \left\{ (\mathbf{m},\mathbf{u}) :\ \mathbf{m} \in L^p(0,T;BV(\Omega)),\ \mathbf{u}\in L^q(0,T;BV(\Omega))^2,\ \nabla \cdot \mathbf{u} \in \Theta \right\},
\end{align*}

\noindent
$\Theta$ being a Lebesgue space including an upper bound constraint. 
\end{theorem}

\begin{proof}
The proof is adapted from the one of \cite[Theorem 3.1]{burger2018} and we just have to prove the weak-$*$ lower semi-continuity property of the additional term $\|\psi \mathbf{m}\|_{L^p(0,T;L^1(\mathbb{R}^4))}$.\\
Since the Wavelet transform is a linear operator in each direction defined by a convolution, then it is a strongly continuous operator from $L^2(\mathbb{R}^2)$ to $L^2(\mathbb{R}^4)$. Also, since the spatial dimension is two, we have the compact embedding of $BV(\Omega)$ into $L^2(\Omega)$ and using the continuous embedding of $L^2(\mathbb{R}^4)$ in $L^1(\mathbb{R}^4)$, 
one can prove the weak-$*$ lower semi-continuity of the additional term $\|\psi \mathbf{m}\|_{L^p(0,T;L^1(\mathbb{R}^4))}$.
\end{proof}

We remark that unlike the models where authors applied directly the $TV-L^1$ optical flow model, we seek to further improve that formulation via a learnt dictionary. To do this, we propose to encode the flow estimation in a sparse representation over a learnt dictionary. Therefore, the discrete formulation of \textit{our proposed joint and hybrid model} after the learning step of the dictionary reads: 

\begin{align}
    \label{eq:actual_model}
    \inf_{\mathbf{m},\mathbf{u},D,\mathbf{a}} E(\mathbf{m},\mathbf{u},D,\mathbf{a}) &= 
    \int_0^T \frac{1}{2}\|K\mathbf{m}-\mathbf{f}\|_{F}^2 + \lambda_1 TV(\mathbf{m}) + \lambda_2\|\psi \mathbf{m}\|_{1} \\ 
    &\quad +\lambda_3\|\frac{\partial \mathbf{m}}{\partial t} + \nabla \mathbf{m} . \mathbf{u} \|_{1} + \lambda_4 TV(\mathbf{u}) \nonumber \\ 
    &\quad + \sum_{\mathbf{p}\in\mathcal{P}}\lambda_5 \|R_{\mathbf{p}}\mathbf{u}-D\mathbf{a}_{\mathbf{p}}\|_F^2+\lambda_6 \|\mathbf{a}_{\mathbf{p}}\|_1\, dt \nonumber
\end{align}

Since all the norms are equivalent in the discrete setting, we have removed the $p$ and $q$ coefficients. We now turn to the numerical resolution of this challenging problem.

%% file: sections/3_algorithm.tex


Due to the non-convexity of the energy coming from the brightness constancy assumption, its non-differentiability with the TV and $L^1$ terms, and complex linear operators applied to the unknowns ($K$ and $\psi$), the problem \cref{eq:actual_model} 
is very challenging to solve numerically. We thus designed an alternating minimisation scheme, to solve our new hybrid model, in which we fix all the variables except one and solve the simplified sub-problem for each unknown in turn. The algorithm is therefore divided into four steps which will be described in the following sub-sections, each one relying on the primal dual Chambolle and Pock algorithm \cite{Chambolle::2011}. We now recall the general setting of this method. The Chambolle and Pock procedure aims at solving the nonlinear primal problem $\min_{x\in X} F(C x)+G(x)$ with the following primal dual formulation:

\begin{align}
    \label{eq:primal-dual}
    \min_{x\in X}\max_{y\in Y}\ \langle C x,y\rangle + G(x) -F^{\star}(y).
\end{align}

The hypothesis are: $X$ and $Y$ are finite dimensional real vector spaces, the map $C:X\rightarrow Y$ is a continuous linear operator, $G:X\rightarrow[0,+\infty[$ and $F^{\star}:Y\rightarrow[0,+\infty[$ are proper, convex, lower-semicontinuous functions and $F^{\star}$ is the Legendre-Fenchel conjugate of $F$. The following \cref{alg:chambolle-pock} converges towards a saddle-point of \cref{eq:primal-dual}.
\begin{algorithm} \caption{Chambolle and Pock iteration} \label{alg:chambolle-pock}
\begin{algorithmic}
    \STATE{Choose $\tau,\sigma>0$ such that $\tau\sigma\|C\|^2<1$, $\theta\in[0,1]$, $(x^0,y^0)\in X\times Y$ and set $\bar{x}^0=x^0$}
    \STATE{Update $x^n,y^n$ and $\bar{x}^n$ as follows:}
    \begin{align}
        \left\{\begin{array}{cl}
            y^{n+1} &= (I+\sigma\partial F^{\star})^{-1}(y^n+\sigma C\bar{x}^n)  \\
            x^{n+1} &= (I+\tau \partial G)^{-1}(x^n-\tau C^{\star}y^{n+1})\\
            \bar{x}^{n+1} &= x^{n+1} + \theta(x^{n+1}-x^n)
        \end{array}\right.
    \end{align}
\end{algorithmic}
\end{algorithm}
The resolvent operators defined through
\begin{align}
    x= (I+\tau\partial F)^{-1}(y) = \arg\min_{x}\left\{\frac{\|x-y\|^2}{2\tau}+F(x)\right\}.
\end{align}
\noindent

Our reconstructed sequence is initialised as the zero filling inverse Fourier reconstruction from the sub-sampled measurements. Zero filling is a common choice for initialisation  as it does not add any substantial computational time. A discussion on other possible initialisations can be seen in Appendix D.
The first approximation of the optical flow $\mathbf{u}^0$ is computed on this image sequence by the $TV-L^1$ method. We consider a sequence of $N_t$ frames with size $N_x\times N_y$ pixels. Each frame will be decomposed in overlapping square patches of size $P_s\times P_s$ and we choose $N_d$, the number of atoms in the dictionary such that $N_d(P_s)^2 > N_xN_y$. We denote by $P_n$ the number of patches for a given frame. Then the main loop of our algorithm is described by \cref{alg:main-loop}.


\begin{algorithm} \caption{Main Loop} \label{alg:main-loop}
\begin{algorithmic}
    \STATE{Given a threshold $\epsilon > 0$.}
    \STATE{Let $\mathbf{m}^0$ be a first approximation of the MRI sequence (ZF).}
    \STATE{Let $\mathbf{u}^0$ be a first approximation of the optical flow ($TV-L^1$).}
    \STATE{Learning the dictionary $D$ (if needed).}
    \STATE{Let $\mathbf{a}^0$ be a first approximation of the sparse decomposition of $\mathbf{u}^0$ in $D$ (from the learning part).}
    \REPEAT
    \STATE{Compute $\mathbf{m}^{n+1}=\arg\min_{\mathbf{m}} E(\mathbf{m}, \mathbf{u}^n, \mathbf{a}^n)$}
    \STATE{Compute $\mathbf{u}^{n+1}=\arg\min_{\mathbf{u}} E(\mathbf{m}^{n+1}, \mathbf{u}, \mathbf{a}^{n})$}
    \STATE{Compute $\mathbf{a}^{n+1}=\arg\min_{\mathbf{a}} E(\mathbf{m}^{n+1}, \mathbf{u}^{n+1}, \mathbf{a})$}
    \UNTIL{$\|\mathbf{m}^{n+1} -\mathbf{m}^n\| < \epsilon\|\mathbf{m}^n\|$}
    \RETURN $\mathbf{m}^{n+1}$.
\end{algorithmic}
\end{algorithm}

In what follows, we describe each of the four required algorithms, namely: learning the dictionary, image reconstruction, optical flow estimation, and sparse coding.

\subsection{Learning the dictionary}
In this section, we present the process of learning the dictionary for the optical flow consisting in the minimisation of \cref{eq:learn_dico} following the strategy of \cite{olshausen1996emergence} which we recall here for the sake of completeness. Given a discrete reference optical flow $\mathbf{u}^{ref}=(u^{ref}_x,u^{ref}_y)^T\in [\mathbb{R}^{N_x\times N_y\times(N_t-1)}]^2$, dictionary learning aims at finding the best dictionary $D_{\mathbf{u}}=\begin{pmatrix}D_x&0\\0&D_y\end{pmatrix}$ with $D_x$ and $D_y$ having $N_d$ atoms such that all the data can be expressed as a linear combination of these atoms. 
Let $R_{\mathbf{p}} = \begin{pmatrix}R_{x,\mathbf{p}}&0\\0&R_{y,\mathbf{p}}\end{pmatrix}$ be the extracting patch operator centred at pixel $\mathbf{p}$. Let $\mathcal{P}$ be the set of all patches, we set $P_n=|\mathcal{P}|$ and $R=(R_{\mathbf{p}})_{\mathbf{p}\in\mathcal{P}}$ so that the extracting patch operator $R$ goes from $[\mathbb{R}^{N_x\times N_y\times(N_t-1)}]^2$ onto $[\mathbb{R}^{P_s^2\times P_n\times(N_t-1)}]^2$. Let $\mathbf{a}=\begin{pmatrix}a_x\\a_y\end{pmatrix} \in [\mathbb{R}^{N_d\times P_n\times(N_t-1)}]^2$ be the sparse coefficients of the decomposition of $\begin{pmatrix}u_x\\u_y\end{pmatrix}$ in $D_{\mathbf{u}}$.

From now on, we keep the notation $u$ for both $u_x$ and $u_y$, the notation $R$ for both $R_x$ and $R_y$, the notation $D$ for both $D_x$ and $D_y$ and the notation $a$ for both $a_x$ and $a_y$ as the learning process is decoupled for both directions and can be performed in parallel computing. 

Let $\bar{u} = Ru^{ref} \in \mathbb{R}^{P_s^2\times P_n\times (N_t-1)}$ be such that $\bar{u}_k \in\mathbb{R}^{P_s^2\times P_n}$ is the matrix of all patch extracted from $\bar{u}$ at time $k$.  And let $a=(a_k)_{k=1}^{N_t-1}$ defined by

\begin{align}
    \min_{\|a\|_{l^0} \leq \epsilon} \sum_{k=1}^{N_t-1} \|\bar{u}_k  - D a_k\|_2^2.
\end{align}

\noindent
The pseudo-norm $l^0$ of $a\in\mathbb{R}^{N_d\times P_n\times (N_t-1)}$ is given by $\|a\|_{l^0}=\underset{1\leq k\leq N_t-1}{\min} |\{(a_k)_{i,j}\neq 0\}|$ and $\epsilon >0$ controls the amount of sparsity. Dictionary learning performs an optimisation both on the dictionary $D$ and on the set of coefficients $a$. This joint optimisation reads

\begin{align}
    \label{eq:learning_energy}
    \min_{D\in \mathcal{D}, a\in\mathcal{A}} E(D,a) &= \sum_{k=1}^{N_t-1} \frac{1}{2}\|\bar{u}_k - Da_k\|_F^2
\end{align}

The constraint set on $D$ is $\mathcal{D}= \{D\in\mathbb{R}^{P_s^2\times N_d},\ \forall\ 1\leq j\leq N_d,\ \|D_j\|\leq 1\}$, the columns of the dictionary are unit normalised, the sparsity constraint set on $a$ reads $\mathcal{A}= \{a\in \mathbb{R}^{D_n\times P_n},\ \forall\ 1\leq j\leq P_n,\ |a_j|_0\leq \epsilon \}$.
\begin{remark}These constraints ensure uniqueness of the solution, indeed if they were not here, we could always find a transformation $P$ such that $P\mathbf{a}$ has smaller $l^0$ norm and the transformed dictionary $DP^{-1}$ has bigger norm. \end{remark}

Following standard arguments, we propose a block coordinate descent method to minimise the energy:

\begin{align}
    D^{n+1} &\in\arg\min_{D\in\mathcal{D}}E(D,a^{n}),\quad a^{n+1} \in\arg\min_{a\in\mathcal{A}}E(D^{n+1},a).
\end{align}

Convergence on such alternating minimisation scheme can be shown, see for instance \cite{tseng2001convergence}.

For the first sub-problem, we assume that $a$ is fixed, then it becomes a quadratic optimisation problem under constraints with respect to $D$:

\begin{align}
    \min_{D\in\mathcal{D}} \left\|\bar{u}_k\left(\sum_{k=1}^{N_t-1} a_k^T\right)\sqrt{\sum_{k=1}^{N_t-1}a_ka_k^T}^{-T} - D \sqrt{\sum_{k=1}^{N_t-1}a_ka_k^T}\right\|_F^2 .
\end{align}

The minimisation problem can be exactly solved as:

\begin{align}
    \label{eq:projection_dictionary}
    D^{n+1} = \text{proj}_{\mathcal{D}}\left(D^n -\tau_1\left(D^n\sum_{k=1}^{N_t-1}a_k a_k^T -\sum_{k=1}^{N_t-1}\bar{u}_ka_k^T  \right) \right)
\end{align}

\noindent
where the projection step is $\tau_1 < \frac{2}{\left\|\sum_{k=1}^{N_t-1}a_ka_k^T\right\|}$ and $\text{proj}_{\mathcal{D}}(D)=\frac{D}{\|D\|}$.

Then we consider $D$ fixed and minimise with respect to $a$. This again leads to a constrained problem which can be solved by projections for each time step on the constraint set,

\begin{align}
    \label{eq:projection_sparse}
    a_{k}^{n+1} &= \text{proj}_{\mathcal{A}}\left(a_{k}^n-\tau_2 D^T (Da_k^n-U_{k})\right) 
\end{align}

\noindent
where this time, $\tau_2 < \frac{2}{\|DD^T\|}$. For $1\leq j\leq P_n$ we denote $|\bar{a}_{k,j}(1)|\leq \cdots\leq |\bar{a}_{k,j}(N_d)|$ the order of magnitude of the vector $a_{k,j}\in\mathbb{R}^{N_d}$. Then the projection operator on $\mathcal{A}$ reads:  

\begin{align}
    \forall 1\leq i\leq N_d,\ \tilde{a}_{k,j}(i) = \left\{\begin{array}{cl}
        a_{k,j}(i) &\text{if}\ |a_{k,j}(i)| \geq |\bar{a}_{k,j}(\epsilon)|  \\
        0 &\text{otherwise.} 
    \end{array}\right.
\end{align}

\begin{remark}
    Two levels of parallel computing can be performed for this process. First we can learn the dictionary $D_x$ for the horizontal component $u$ of the optical flow and the dictionary $D_y$ for the vertical component simultaneously. Additionally, for each component we can minimise $a_k$ in parallel. 
\end{remark}

\subsection{Image reconstruction} 
This section describes the process of MRI reconstruction which is the core of the method. While in \cite{lustig2007} a smooth regularisation of the $L^1$ norm was used, in this work, we propose to solve the non-smooth optimisation problem using \cref{alg:chambolle-pock}. Assume $\mathbf{u}$ is known, the problem then reads

\begin{align}
    \label{eq:image_reconstruction}
    \min_{\mathbf{m}} E(\mathbf{m}) = \int_0^T \frac{1}{2}\left\|K\mathbf{m} - \mathbf{f}\right\|^2 + \lambda_1\left\|\nabla\mathbf{m}\right\|_{2,1} + \lambda_2\left\|\Psi\mathbf{m}\right\|_1 +\lambda_3\left\|\frac{\partial\mathbf{m}}{\partial t} + \nabla\mathbf{m}\cdot\mathbf{u}\right\|_1\ dt
\end{align}

\noindent
with $\|\nabla \cdot\|_{2,1}$ representing the discrete spatial isotropic $TV$. 

We proceed as in \cite{burger2018}. Let $C=\begin{bmatrix} K,\ \nabla,\ \Psi,\ D_t + \mathbf{u}\cdot\nabla \end{bmatrix}^T$ be the linear operator acting on $\mathbf{m}$ so that the minimisation is seen as $\min_{\mathbf{m}}E(C\mathbf{m})$. And let $\mathbf{y}=\begin{bmatrix} y_0,y_1,y_2,y_3 \end{bmatrix}^T$ be the collection of dual variables for each operator. From \cref{eq:primal-dual} we have $F:=E$ and $G=0$ and the Legendre-Fenchel conjugate of the energy is given by 

\begin{align}
    E^{\star}(\mathbf{y}) = \int_0^T &\frac{1}{2}\|y_0\|_2^2+\langle y_0,\mathbf{f}\rangle + \delta_{\{y:\|y\|_{2,\infty}\leq 1\}}\left(\frac{y_1}{\lambda_1}\right) + \delta_{\{y:\|y\|_{\infty}\leq 1\}}\left(\frac{y_2}{\lambda_2}\right)\\ 
\nonumber &+\delta_{\{y:\|y\|_{\infty}\leq 1\}}\left(\frac{y_3}{\lambda_3}\right) dt
\end{align}

\noindent
with $\delta_I(y) = 0$ if $y\in I $ and $+\infty$ otherwise. The primal-dual algorithm requires the proximity operator $(I+\sigma\partial E^{\star})^{-1}(\mathbf{y}^{n}+\sigma C\bar{\mathbf{m}}^{n})$. The components $y_i^{n+1}$ for $i=0,1,3$ are computed exactly as in \cite{burger2018} and we give only 

\begin{align*}
    y_2^{n+1} &= \arg\min_{y}\left\{\int_0^T\frac{1}{2}\left\|y- y_2^n-\sigma\Psi\left(\mathbf{m}^{n} + \theta(\mathbf{m}^{n} -\mathbf{m}^{n-1})\right) \right\|_2^2 +\delta_{\{y:\|y\|_{\infty}\leq 1\}}\left(\frac{y}{\lambda_2}\right) dt\right\}.
\end{align*}

This subproblem admits a closed form solution recalled here,

\begin{align*}
    \tilde{\mathbf{y}}^{n+1} &= \mathbf{y}^{n}+\sigma C\left(\mathbf{m}^{n} + \theta(\mathbf{m}^{n} -\mathbf{m}^{n-1} )\right) \\ 
    \mathbf{y}^{n+1} &= (I + \sigma\partial E^{\star})^{-1}(\mathbf{y}^{n}) = 
    \left\{ 
        \begin{array}{ccl}
        y_0^{n+1} &= \frac{y_0^{n} - \sigma \mathbf{f}}{\sigma +1}, \\ 
        y_1^{n+1} &=\pi_{\lambda_1}(\tilde{y}_1^{n} ),\\ 
        y_2^{n+1} &=\pi_{\lambda_2}(\tilde{y}_2^{n} ),\\
        y_3^{n+1} &=\pi_{\lambda_3}(\tilde{y}_3^{n} )
        \end{array} .
    \right.
\end{align*}

The projection onto the unit ball is given by $\pi_{\lambda}(y) = \frac{y}{\max\left(1,\frac{\|y\|_2}{\lambda}\right)}$. The parameters $\tau$, $\sigma$ refer to step sizes, and have to fulfil the following constraints: $\tau\sigma\|C\|^2\leq 1$ with $\|C\| \approx 2+\sqrt{8}+\sqrt{2}(1+\|\mathbf{u}\|)$.

\subsection{Sparse flow representation and optical flow approximation}

The minimisation over the variable $\mathbf{a}$ can be split for each component $\mathbf{a}_{\mathbf{p}}$ which is given by $\tilde{E}(\mathbf{a}_{\mathbf{p}}) = \lambda_5\|R_{\mathbf{p}}\mathbf{u}- D\mathbf{a}_{\mathbf{p}}\|_F^2 + \lambda_6\|\mathbf{a}_{\mathbf{p}}\|_1$. From \cref{alg:chambolle-pock} we take $F=\lambda_6\|\cdot\|_1$, the proximal operator of its Legendre-Fenchel conjugate is given by projection operator: 

\begin{align}
    \label{eq:prox_sparse_alg}
    (I+\sigma \partial F^{\star})^{-1} (\mathbf{y}) = \pi_{\lambda_6} (\mathbf{y})
\end{align}

Moreover since $G=\lambda_4\|R_{\mathbf{p}}\mathbf{u}- D\cdot\|_2^2$ is smooth, its resolvent reduces to the gradient, 

\begin{align}
    \label{eq:grad_prox_alg}
    \tau \nabla G(\mathbf{a}_{\mathbf{p}}) = 2\tau\lambda_5 D^T(D\mathbf{a}_{\mathbf{p}} - R_{\mathbf{p}}\mathbf{u}).
\end{align}

Finally we optimise on the optical flow $\mathbf{u}$, given fixed $\mathbf{m}$ and $\mathbf{a}=(\mathbf{a}_{\mathbf{p}})_{\mathbf{p}\in\mathcal{P}}$. Let $E(\mathbf{u}) = F(\nabla \mathbf{u}) + G(\mathbf{u})$ with $F(\cdot)=\lambda_4\|\cdot\|_{2,1}$ and $G=\lambda_3\left\|\frac{\partial \mathbf{m}}{\partial t}+\nabla\mathbf{m}\cdot \right\|_1+\lambda_5\sum_{\mathbf{p}}\|R_{\mathbf{p}}\cdot-D\mathbf{a}_{\mathbf{p}}\|_F^2$. For simplicity of the notations we introduce the operators $\mathcal{A}=I+2\tau\lambda_5\sum R_{\mathbf{p}}^T R_{\mathbf{p}}$,  $\rho(\mathbf{u}) = \frac{\partial \mathbf{m}}{\partial t}+\mathbf{u}\cdot\nabla\mathbf{m}$ and the notation $\tilde{\mathbf{u}}=\mathbf{u}+2\tau\lambda_5\sum R_{\mathbf{p}}^TD\mathbf{a}$. Then the proximal operator of $G$ reads, 

\begin{align}
    \label{eq:optflow_prox}
    \mathcal{A}(I+\tau\partial G)^{-1} (\mathbf{u}) &= \tilde{\mathbf{u}} + 
    \left\{
        \begin{array}{cll}
            -\lambda_3 \nabla\mathbf{m} & \text{if} & \rho(\mathcal{A}^{-1}\tilde{\mathbf{u}}) > \tau\lambda_3 \mathcal{A}^{-1}\|\nabla\mathbf{m}\|^2, \\
            \lambda_3\nabla\mathbf{m}   & \text{if} & \rho(\mathcal{A}^{-1}\tilde{\mathbf{u}}) < - \tau\lambda_3 \mathcal{A}^{-1}\|\nabla\mathbf{m}\|^2, \\
             -\rho(\tilde{\mathbf{u}}) \frac{\nabla\mathbf{m}}{\|\nabla\mathbf{m}\|^2}& \text{else,}
        \end{array}
    \right.
\end{align}

and the $F^{\star}$ one is

\begin{align}
    \label{eq:optflow_prox_2}
    (I+\sigma\partial F^{\star})^{-1}(\mathbf{y}) &= \pi_{\lambda_6}(\mathbf{y}).
\end{align}

Details of the computation of \cref{eq:optflow_prox} can be found in \cref{anx:proximal} and details about the discretisation of the operators can be found in \cref{anx:operators}. We now turn to the numerical experiments. 

%% file: sections/4_numericals.tex
In this section we show the results of our numerical experiments. All tests and comparisons were run under the same CPU-based implementation (Intel Core i7-8550U) with MATLAB 2018b. Compared techniques were either used as provided by the author or implemented when not available. 

\subsection{Data Description and Evaluation Protocol}\label{sub-sec:data_description}

\begin{figure*}[!t]
\centering
\includegraphics[width=1\textwidth]{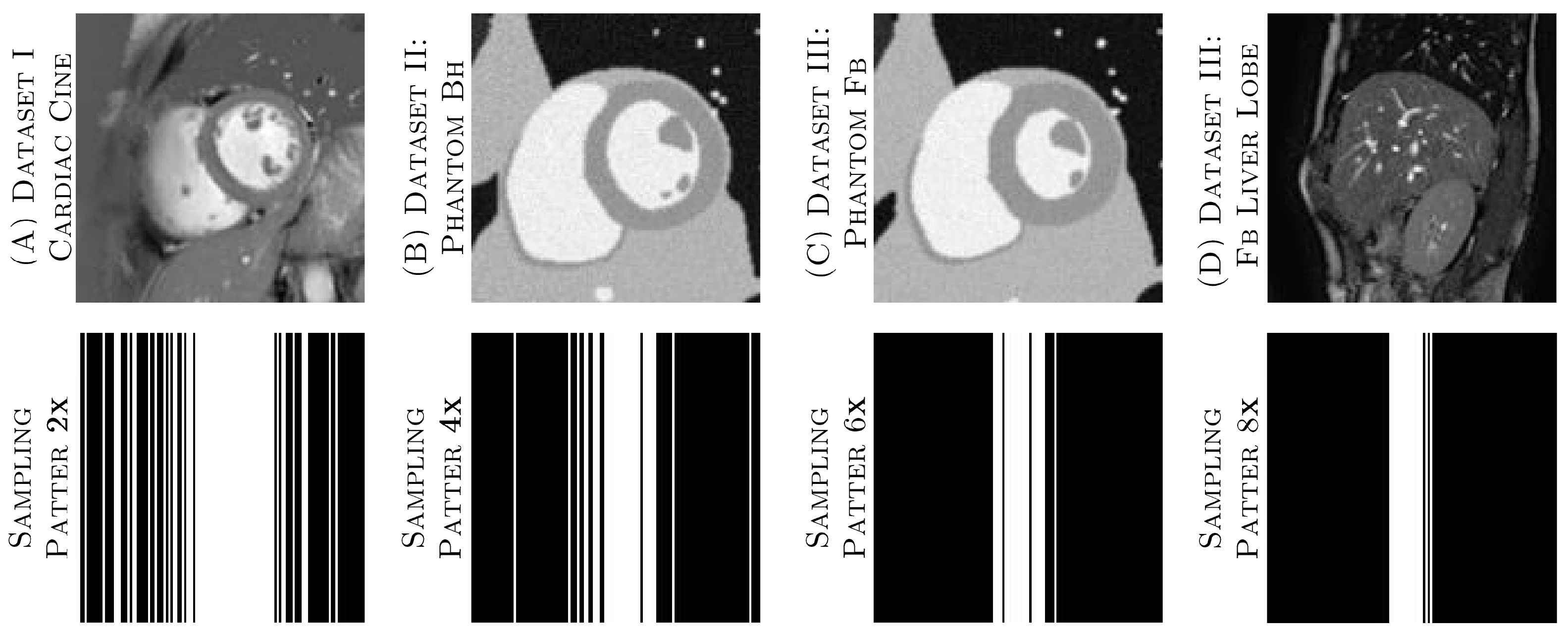}
\caption{(Top column) Visual samples of the datasets used in our experiments. (Bottom row) visualisation of some undersampling patterns used in our experiments using acceleration factor=\{2x,4x,6x,8x\}.} 
\label{fig::visu1}
\end{figure*}

Our approach is evaluated by using data coming from: cine cardiac, right lobe liver and phantom generated cardiac cine. The datasets were saved as fully sampled raw data and then were retrospectively undersampled using a variable-density random sampling pattern, and in particular, a cartesian pattern. \footnote{Raw data and method to create the undersampled datasets can be download here \url{https://github.com/tschmoderer/2019_mri_reconstruction/blob/master/src/data/data_mri.tar.xz}} \cref{fig::visu1} as suggested by Lustig \cite{lustig2007}. Datasets characteristics are as follows:

\begin{itemize}
    \item Dataset I: 
    A cine cardiac dataset acquired from a healthy volunteer, from \cite{Schloegl::2017}. Matrix size - $128\times128$ pixels and 30 frames.
    \item Dataset II and III: 
    Realistic cardiac cine simulation generated using the MRXCAT phantom framework \cite{Wissmann::2014}. Whilst the Dataset II was simulated during breath holding, the III was set with free respiratory motion. Matrix size - $128\times128$ and 24 frames. 
    \item Dataset IV: 
    4DMRI dataset~\cite{von20074d} acquired from the right liver lobe during free-brea-thing. It was acquired using a 1.5T Philips Achieva system with: TR=3.1ms,  slices=25 and matrix size of $195\times166$.  
\end{itemize}

In this work, we assumed single coil datasets. However, our work can be extended to the multiple coils case, which will be the focus of our future work.

We have followed common protocol for evaluating our MRIR-DLMC using  two metrics over the full reconstructions~\cite{lustig2007,Lingala::2015,LustigktSPARSE::2006,Majumdar::2015}: the peak signal-to-noise ratio (PSNR) and the structural similarity (SSIM). During the following sections we use MRIR-DLMC and OURS to refer to our technique. The metric-wise computations were performed between the gold-standard (i.e. fully sampled case) and the reconstructed MR images. The explicit definition of these metrics can be found in \cref{anx:err_measures}. To demonstrate the potentials of our approach, we designed a two-part evaluation protocol. The first one consists in comparing our technique versus three top reference techniques for MRI reconstructions: i) \emph{zero-filling} (ZF) approach i.e. an inverse Fourier transform on the undersampled data ii) the classic \emph{compressed sensing} (CS) algorithm \cite{lustig2007} and iii) the \emph{low rank and sparsity} (L+S) technique of~\cite{Otazo::2015}. The second comparison scheme is based on techniques that perform also two tasks simultaneously: \emph{compressed sensing and motion} (CS+M) \cite{aviles2018} and the recent method of (JPDAL) \cite{Zhao::2019}. The parameters were set based on an empirical testing in a coarse to fine search strategy by maximising the SSIM metric. From this, we set the parameter for our experiments as displayed in Table~\ref{tab:parameters}.


\begin{table}[t!]
    \centering
    \resizebox{\textwidth}{!}{  
    \begin{tabular}{@{}ccccc@{}}
    \toprule
    \multicolumn{1}{c}{\begin{tabular}[c]{@{}c@{}} \textsc{Weighing} \\\textsc{Parameters} \cref{eq:actual_model}\end{tabular}}&
      \multicolumn{1}{c}{\begin{tabular}[c]{@{}c@{}}\textsc{Image} \\  \textsc{reconstruction} \\ (\cref{alg:chambolle-pock} for \cref{eq:image_reconstruction})\end{tabular}} &
      \multicolumn{1}{c}{\begin{tabular}[c]{@{}c@{}}\textsc{Sparse} \\  \textsc{Approximation} \\ (\cref{alg:chambolle-pock} with \cref{eq:prox_sparse_alg,eq:grad_prox_alg})\end{tabular}} &
      \multicolumn{1}{c}{\begin{tabular}[c]{@{}c@{}}\textsc{Optical Flow} \\  \textsc{Approximation} \\ \cref{eq:optflow_prox,eq:optflow_prox_2}\end{tabular}} &
      \multicolumn{1}{c}{\begin{tabular}[c]{@{}c@{}}\textsc{Dictionary} \\  \textsc{learning} \\ \cref{eq:learning_energy}\end{tabular}} \\ \midrule
    \begin{tabular}[c]{@{}l@{}}$\lambda_1=0.003$\\ $\lambda_2=0.0001$\\ $\lambda_3=0.001$\\ $\lambda_4=0.001$\\ $\lambda_5=0.001$\\ $\lambda_6=0.0001$\\ $\epsilon = 0.001$\end{tabular} &
      \begin{tabular}[c]{@{}l@{}}$\sigma =0.05$\\ $\tau = 0.05$\\ $\theta = 1$\end{tabular} &
      \begin{tabular}[c]{@{}l@{}}$\sigma =0.99$\\ $\tau = 0.99$\\ $\theta = 1$\end{tabular} &
      \begin{tabular}[c]{@{}l@{}}$\sigma =0.5$\\ $\tau = 0.25$\\ $\theta = 1$\end{tabular} &
      \begin{tabular}[c]{@{}l@{}}$\epsilon = 0.7N_d$\\ $\tau_1$ like in \cref{eq:projection_dictionary}\\ $\tau_2$ like in \cref{eq:projection_sparse}\end{tabular} \\ \hline
    \end{tabular}
    }
    \caption{Parameter values used in our experiments.}
    \label{tab:parameters}
\end{table} 

\begin{figure*}[!t]
\centering
\includegraphics[width=1\textwidth]{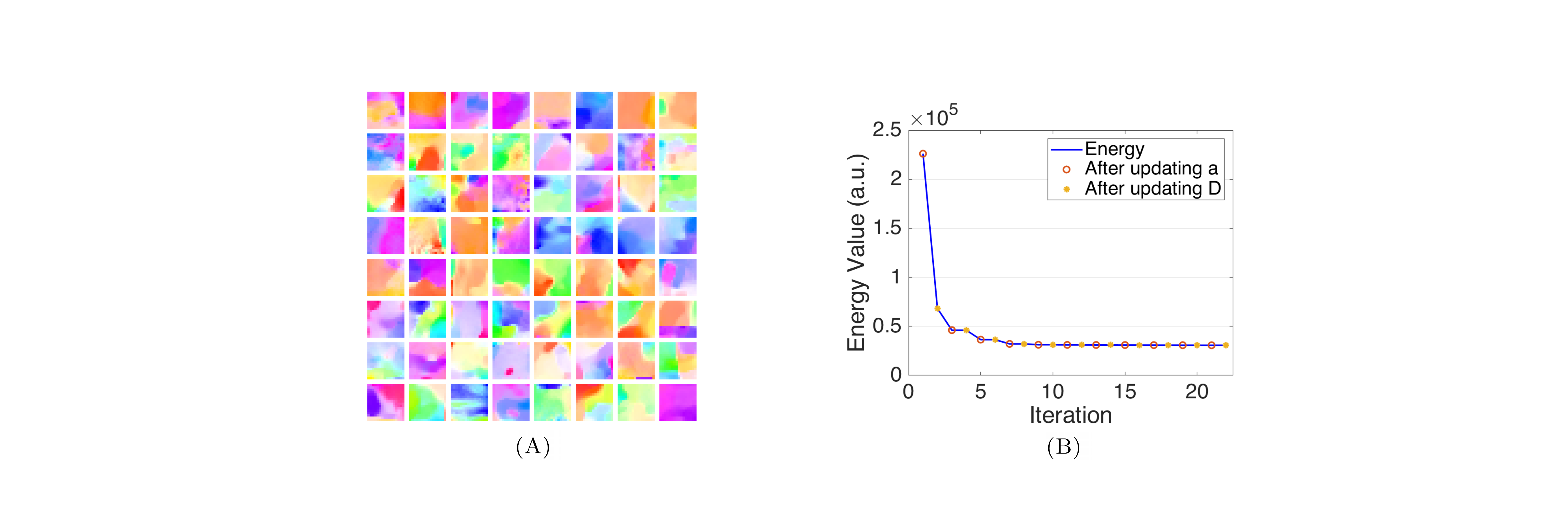}
\caption{(A)  Samples extracted from our learned dictionary with $1024$ atoms and patches of size $16\times16$. (B) Evolution of the energy during the learning process.}
\label{fig::dic}
\end{figure*}

For each dataset, several experiments are performed depending on the utilisation of the dictionary. First, for each dataset and acceleration factor, we train a dictionary based on a $TV-L^1$ approximation of the optical flow computed on the inverse Fourier transform of undersampled data. In a second step, to show the transfer learning capability of our model, we train the dictionary on a $TV-L^1$ approximation of the optical flow from fully sampled phantom data, and use it for real cine cardiac applications. 

\textit{We remark that one of the strengths of our approach is to avoid the need to have a ground truth to construct the dictionary}. This is a highly desirable property 
as in the medical domain, the assumption of having ground truth is strong. We use over-complete patch-based dictionary using the $TV-L^1$ OF approximation. 
The dictionary training process involves the computation of the sparse decomposition of the reference flow, so we use this procedure to initialise our sequence $\mathbf{a}^0$.  
Our dictionary is composed of $N_d=1024$ atoms, but for visualisation purposes we display $64$ atoms illustrated on the left of \cref{fig::dic}. The samples are coloured according to the standard visualisation of optical flow. The plot on the right of \cref{fig::dic} displays the evolution of the energy through the learning process to illustrate the speed and stability of convergence. We remark a fast convergence towards a dictionary $D$ and a sparse decomposition $\mathbf{a}$. This translates into a small computational time required to get the optimal dictionary and sparse decomposition. 



\subsection{Results \& Discussion} \label{sub-sec:numerical results}

\begin{figure*}[hbtp]
\centering
\includegraphics[width=1\textwidth]{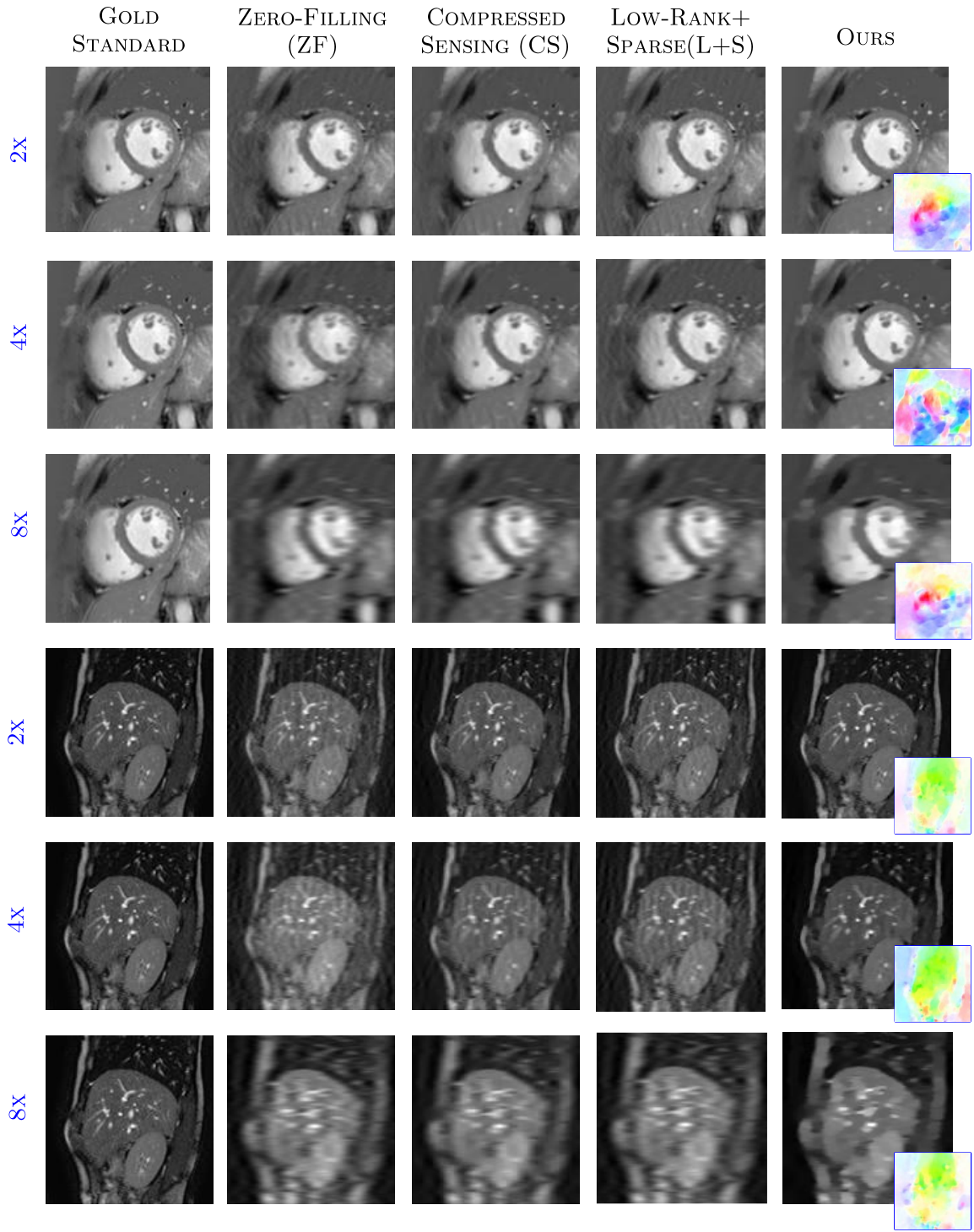}
\caption{Visual comparison of image quality of our model (OURS) against state of the art reconstruction techniques namely: zero filling (ZF), compressed sensing \cite{lustig2007} (CS), and low rank and sparse decomposition \cite{Otazo::2015}, along with the gold standard fully-sampled reconstruction for two real datasets I and IV. For each dataset, the results are shown for three acceleration factors 2x, 4x and 8x. }
\label{fig::res2}
\end{figure*}

\begin{figure*}[hbtp]
\centering
\includegraphics[width=1\textwidth]{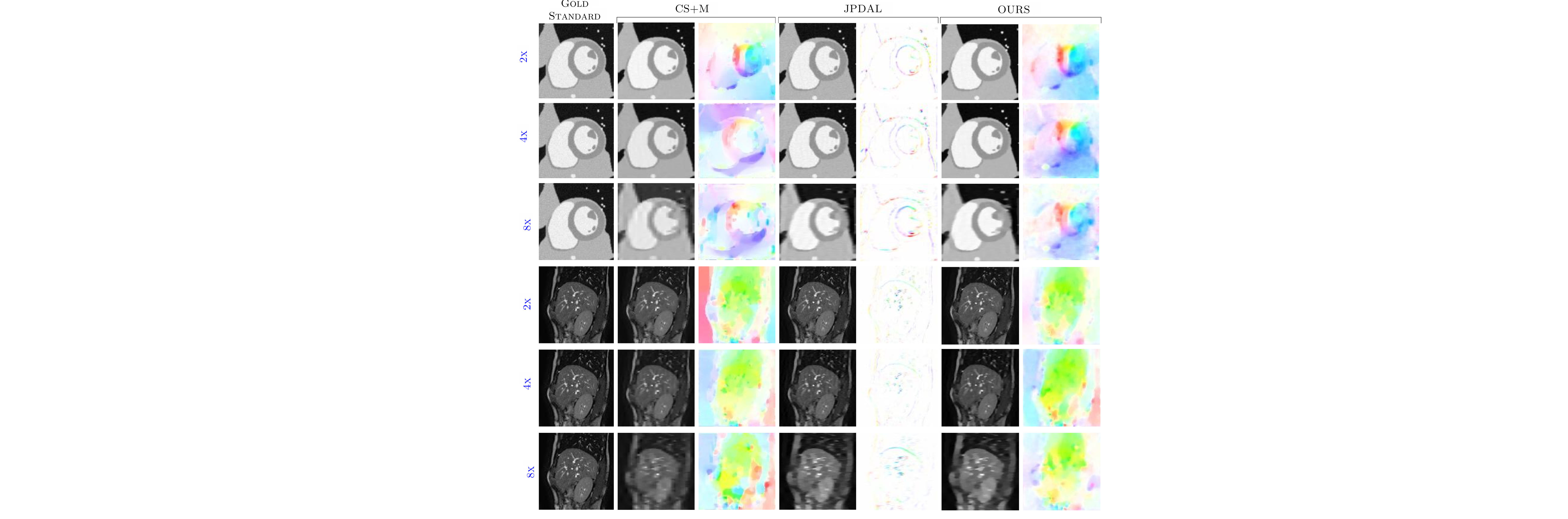}
\caption{Visual comparison of both reconstruction and optical flow estimation quality of our model (OURS) against two joint models namely the purely variational compressed sensing and optical flow model \cite{burger2018,Aviles-Rivero::2018} (CS+M) and the technique introduced in \cite{Zhao::2019} (JPDAL), along with the gold standard fully-sampled reconstruction for datasets II (cardiac phantom) and IV (liver MRI). For each dataset, the results are shown for three acceleration factors 2x, 4x and 8x.} 
\label{fig::res5}
\end{figure*}

\medskip
In this section, we extensively evaluate our approach following the evaluation protocol stated in Section 4.1. 

\textbf{Comparison against single task MRI reconstruction methods.}
In \cref{fig::res2} and following common protocol for MRI reconstruction evaluation, we report the reconstructed images obtained for two real datasets I and IV and three acceleration factors 2x, 4x, and 8x, for a visual comparison of our model against state of the art single task MRI reconstruction techniques: zero-filling (ZF), compressed sensing \cite{lustig2007} (CS) and low-rank and sparse decomposition \cite{Gao::2012} (L+S). We remark that the comparison between our technique and no motion can be observed when comparing the results from the CS reconstructions against ours. Also, the impact of the dictionary learning term can be measured by comparing the CS+M results and ours.
As expected, the zero-filling method performs the worst in all cases and exhibit blurring and wave-like artefacts. We remark that our method outperforms the CS and L+S techniques at all acceleration factors. For an acceleration factor of 2x, this is particularly visible in the region beneath the ventricular chambers of the heart in which both CS and L+S reconstructions still exhibit wave-like and blurring artefacts contrary to ours. As the acceleration factor increases, the difference is even more visible as our method reduces the blurring and wave-like artefacts, better preserves fine details and exhibits sharper edges. Especially, the white blood vessels of the liver are better recovered with our technique than with the other presented ones. The edges of the right and left ventricular chambers of the heart as well as of the liver are sharper with our algorithm than with purely reconstruction methods specially for an acceleration factor of 8x. The blurring and wave-like artefacts in the liver and the myocardium are reduced with our technique for acceleration factors 4x and 8x compared to other reconstruction methods. All these statements are in favor of incorporating motion information in the reconstruction process of dynamic sequences to improve the overall image quality.

\medskip
\textbf{Comparison against joint MRI reconstruction and motion estimation models.}
We now compare our method with other existing joint models addressing simultaneously MRI reconstruction and motion estimation namely a purely variational compressed sensing and optical flow model \cite{burger2018,Aviles-Rivero::2018}, named CS+M, to evidence the improvement induced by the dictionary learning process, and the new compressed sensing and affine optical flow model introduced in \cite{Zhao::2019} called JPDAL. Both the reconstruction and the motion estimation are reported in \cref{fig::res5} for datasets II and IV and for three acceleration factors 2x, 4x and 8x, along with the gold standard fully-sampled reconstruction. Regarding the motion maps, we observe that the JPDAL retrieves small displacements mainly located on the boundary of the organs whereas our model is able to estimate bigger displacements located in the whole organ which seems to be closer to the physiological motion involved. Also, by visual comparison, adding the dictionary learning process tends to smooth and to better identify small moving regions compared to the purely variational model. The motion estimation improvement is reflected in higher quality and visually more pleasing reconstructed images especially for large acceleration factors. For acceleration factors 2x and 4x, the reconstruction for all the methods are visually comparable and competitive. For an acceleration factor of 8x, the CS+M reconstructions exhibit blurring artefacts especially in the liver and the left ventricular chamber of the heart which are less visible in our results. This is probably due to residual movements and less accurate motion estimation and shows the benefits of adding the dictionary learning process in the model to improve image quality. Small details such as the white blood vessels are also better recovered with our technique than with the CS+M one. The JPDAL reconstructions for an acceleration factor of 8x are also more blurred than ours and the homogeneous regions of the left ventricular chamber and the above area are better recovered with our algorithm than with the JPDAL model.

\begin{figure*}[hbtp]
\centering
\includegraphics[width=1\textwidth]{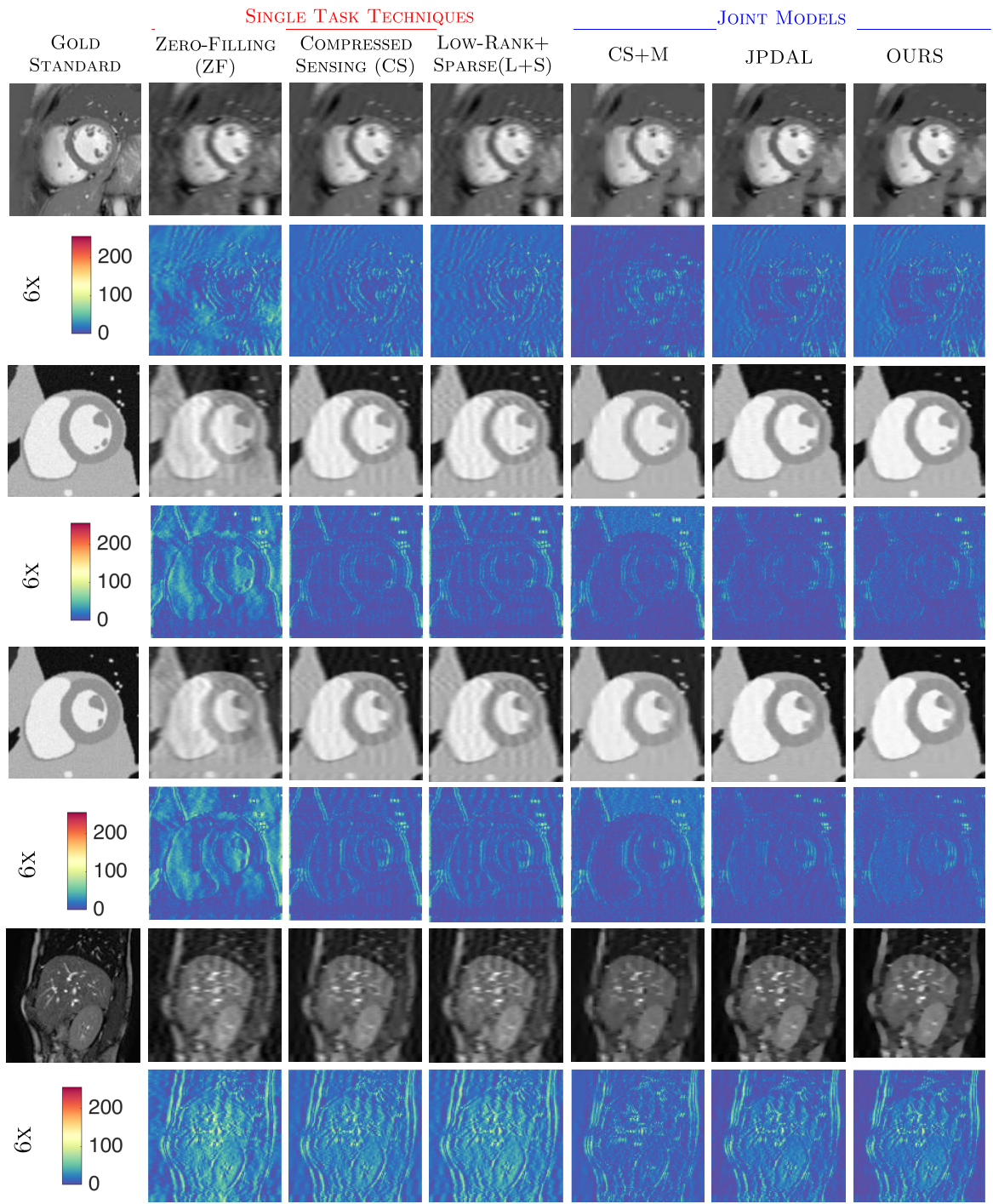}
\caption{Visual comparison of our technique versus single and other joint models. The results are displayed using 6x acceleration factor. A closer look at the difference maps, one can see that our technique reports reconstructions closer to the gold-standard with respect to single task techniques whilst readily compete with, and several times outperforms, the other compared joint models. } 
\label{fig::res1.pdf}
\end{figure*}

\medskip
\textbf{Global Evaluation Performance.} We begin by reporting visual outputs of our proposed algorithmic technique versus single and other joint models along with the gold-standard (fully-sampled case). The results are displayed in  
\cref{fig::res1.pdf} using a high undersampling factor 6x.
By visual inspection, we observe that the outputs generated by the single task techniques (i.e. ZF, CS and LS) are the ones displaying more artefacts than ours, and in particular, as expected, the worst performance is reported by ZF. This in terms of blurring artefacts and loss of preservation of relevant anatomy parts - for example in the cardiac datasets this can be observed around the papillary muscle of the heart whilst in the last dataset in the interior of the liver lobe. These observations are further supported by the corresponding difference maps, in which our technique reports closer reconstruction to the gold-standard than the compared single task techniques. For further support of our proposed technique, we also display selected outputs of our approach versus other techniques with similar philosophy as ours (i.e. joint models). A closer look at these outputs we can see that while CS+M and JDPAL perform better than the single task techniques, our proposed approach readily competes with these techniques, and in several cases, outperforms them. To further support our technique, we ran a nonparametric statistical test (Wilcoxon test) for two group comparison. From the test we found that our technique is statistically  significantly different than other techniques including JPDAL (significance level = 0.05; p-value= 0.024). Moreover, we also reported a comparison in terms of computational time (in seconds), the results are displayed in Figure 6. These plots compare the computational time for all datasets and different undersampling factors for our technique and other multi-task models. From an inspection to these plots, one can see that our method not only improves performance wise but also demands lower computational time that other existing techniques following similar philosophy.

\begin{figure*}[hbtp]
\centering
\includegraphics[width=1\textwidth]{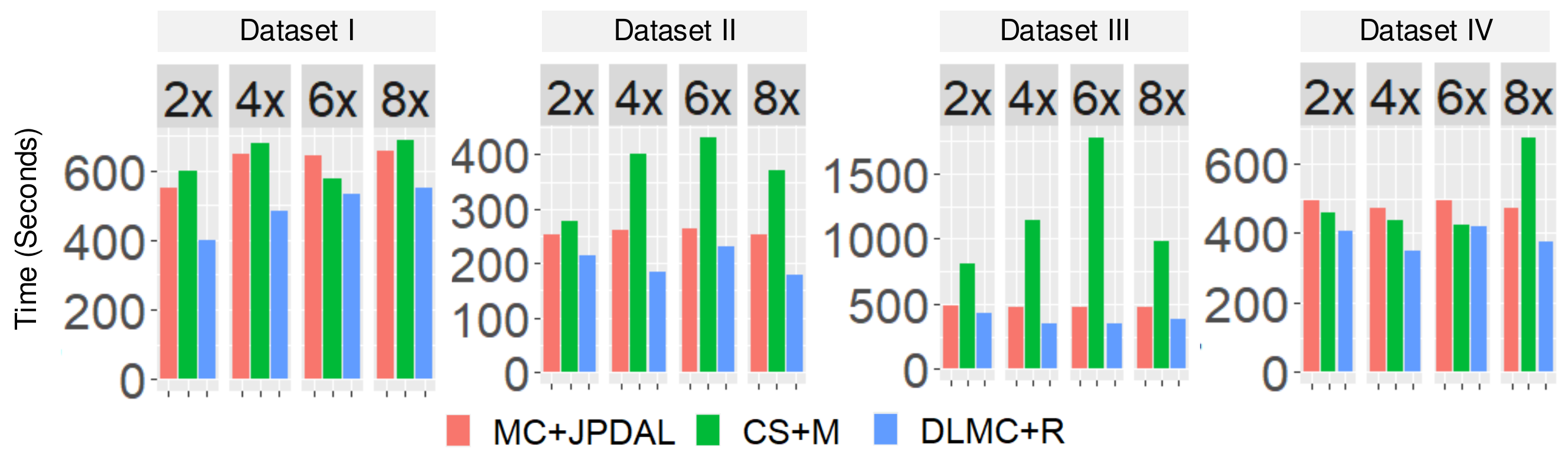}
\caption{Computational time comparison of our approach versus other multi-task techniques reported in seconds. The results are reported from several acceleration factors (2x, 4x, 6x and 8x) and for all datasets. }
\label{fig::comutationalTime}
\end{figure*}
\input{img/results/result_table}

For a more detailed quantitative analysis,  we report the global results in Table~\ref{tab:numerical_results} for all compared techniques and for acceleration factor of 2x, 4x, 6x and 8x. The displayed results are the average  of the image metrics, i.e. PSNR and SSIM, across the full corresponding dataset. With respect to the realistic datasets, we reported the overall best performance whilst for the simulated datasets our method performs similarly  well, and places second only behind JPDAL. We further support our technique by computing the computational time (in seconds) of our technique and other multi-task techniques: CS+M and MC+JPDAL. The results are reported in Figure~\ref{fig::comutationalTime}. In a closer inspection at this figure, one can observe that our technique required for all cases less computational resources than the compared techniques.

\input{img/results/phantomfb_results_table.tex}
\subsection{Transfer Learning Capabilities of our Technique}
We have shown in the previous section the advantages of our technique in terms of quality and computational performance. Besides these benefits, robustness of our training process when no ground-truth is available has been evidenced which is a highly desirable property. In this section, we highlight another strength of our technique, which is the ability of learning over phantom datasets and transfer this knowledge to real datasets. The ability to transfer knowledge from phantom to real data is very challenging yet highly desirable for any algorithmic technique. As pointed out in~\cite{raghu2019transfusion}, this capability is not straightforward and such generalisation might be difficult to achieve due to the overparametrisation of the model given by the fundamental gap between synthetic and real data. 
To prove this capability, we show a numerical comparison in \cref{tab:phantom_fb2cine_numerical_results} and denote the new results as 'DLMCR w/TL'. From the results, one can see that firstly, our approach is capable to perform similarly well, when trained with phantoms datasets, than JPDAL and ours (trained with realistic dataset). Our 'DLMCR w/TL' readily compete numerically with the best scores and in some cases outperforms the compared techniques. This capability opens several opportunities in the medical domain as we avoid the need to retrain our dictionary. 

Overall, our proposed model offers three major advantages over the compared techniques  i) it improves performance wise over the compared techniques, and in particular, over JPDAL with a statistical difference, ii) it requires less computational time than the other techniques and iii) our technique can narrow the gap between synthetic and real data- as our technique performs well when trained on phantoms and then applied to real data.

%% file: img/results/result_table.tex

\begin{table}[]
\resizebox{\textwidth}{!}{  
\begin{tabular}{@{}cccccccccc@{}}
\toprule
\multirow{2}{*}{\textsc{Dataset}} &
  \multirow{2}{*}{\begin{tabular}[c]{@{}c@{}}\textsc{Reconstruction}\\ \textsc{Scheme}\end{tabular}} &
  \multicolumn{2}{c}{2x} &
  \multicolumn{2}{c}{4x} &
  \multicolumn{2}{c}{6x} &
  \multicolumn{2}{c}{8x} \\ \cmidrule(l){3-10} 
 &
   &
  PSNR &
  SSIM &
  PSNR &
  SSIM &
  PSNR &
  SSIM &
  PSNR &
  SSIM \\ \midrule
\multirow{6}{*}{Dataset I}
&
  Zero-Filling &
  31.40 &
  90.19 &
  25.14 &
  81.40 &
  25 &
  76.97 &
  22.99 &
  72.71 \\
 &
  CS &
  32.85 &
  93.58 &
  31.57 &
  88.57 &
  27.98 &
  81.48 &
  22.97 &
  72.76 \\
 &
  L+S &
  34.21 &
  92.77 &
  31.09 &
  86.20 &
  27.74 &
  80.10 &
  22.98 &
  72.68 \\
 &
  CS+M &
  36.72 &
  96.23 &
  31.53 &
  90.26 &
  \textbf{28.39} &
  80.17 &
  \textbf{24.9} &
  72.68 \\
 &
  MC+JPDAL &
  36.7 &
  \textbf{97.85} &
  \textbf{32.72} &
  92.06 &
  27.80 &
  84.29 &
  23.15 &
  75.51 \\
 &
  MRIR-DLMC &
 \textbf{38.01} &
  97.33 &
  32.35 &
  \textbf{92.26} &
  27.65 &
  \textbf{84.76} &
  23.12 &
 \textbf{76.03} \\ \midrule
\multirow{6}{*}{Dataset IV}&
  Zero-Filling &
  22.43 &
  72.95 &
  17.84 &
  57.66 &
  18.49 &
  50.80 &
  16.84 &
  44.98 \\
 &
  CS &
  26.61 &
  83.48 &
  22.7 &
  69.31 &
  19.36 &
  54.15 &
  17.03 &
  45.86 \\
 &
  L+S &
  24.26 &
  77.68 &
  21.71 &
  63.39 &
  19.19 &
  52.46 &
  16.82 &
  44.96 \\
 &
  CS+M &
  31.91 &
  91.73 &
  25.86 &
  76.32 &
  21.14 &
  58.37 &
  20.01 &
  49.17 \\
 &
  MC+JPDAL &
  32.81 &
  93.07 &
  27.3 &
  82.01 &
  22.34 &
  65.62 &
  19.39 &
  54.1 \\
 &
  MRIR-DLMC &
  \textbf{34.21} &
  \textbf{94.16} &
  \textbf{28.28} &
  \textbf{84.26} &
  \textbf{22.65} &
  \textbf{67.17} &
  \textbf{19.57} &
  \textbf{55.02} \\ \midrule
\multirow{6}{*}{Dataset II}&
  Zero-Filling &
  28.27 &
  78.87 &
  23.97 &
  68.08 &
  20.53 &
  61.39 &
  22.57 &
  61.09 \\
 &
  CS &
  31.91 &
  83.11 &
  28.58 &
  76.35 &
  24.83 &
  67.01 &
  22.58 &
  61.12 \\
 &
  L+S &
  30.59 &
  82.17 &
  27.74 &
  73.7 &
  24.63 &
  65.35 &
  22.58 &
  60.91 \\
 &
  CS+M &
  32.65 &
  87.1 &
  30.5 &
  81.91 &
  24.84 &
  70.30 &
  21.94 &
  63.06 \\
 &
  MC+JPDAL &
  \textbf{36.5} &
  \textbf{93} &
  \textbf{32.45} &
  \textbf{86.44} &
  \textbf{28.08} &
  \textbf{76.92} &
  24.46 &
  69.45 \\
 &
  MRIR-DLMC &
  33.19 &
  88.46 &
  30.5 &
  82.77 &
  27.76 &
  76.18 &
  \textbf{25.14} &
  \textbf{71.9} \\ \midrule
\multirow{6}{*}{Dataset III}&
  Zero-Filling &
  28.84 &
  79.3 &
  24.31 &
  67.84 &
  20.57 &
  60.48 &
  22.46 &
  59.61 \\
 &
  CS &
  31.19 &
  82.51 &
  28.66 &
  74.89 &
  25.26 &
  65.87 &
  22.48 &
  59.61 \\
 &
  L+S &
  30.59 &
  82.8 &
  27.99 &
  72.93 &
  24.9 &
  64.5 &
  22.46 &
  59.38 \\
 &
  CS+M &
  32.11 &
  86.13 &
  30.02 &
  \textbf{89.16} &
  24.28 &
  67.74 &
  21.55 &
  60.59 \\
 &
  MC+JPDAL &
  \textbf{35.37} &
  \textbf{91.86} &
  \textbf{31.24} &
  83.85 &
  \textbf{27.3} &
  \textbf{74.16} &
  \textbf{24.14} &
  \textbf{66.66} \\
 &
  MRIR-DLMC &
  32.76 &
  87.48 &
  29.41 &
  80.12 &
  26.9 &
  73.43 &
  24.05 &
  66.5 \\ \bottomrule
\end{tabular}
}
\caption{Numerical comparison of our technique vs single and joint technique for different acceleration factors. The results are reported as the average of the corresponding metric over all the corresponding dataset. }
\label{tab:numerical_results}
\end{table}

%% file: img/results/phantomfb_results_table.tex
\begin{table}[]
 \resizebox{\textwidth}{!}{  
\begin{tabular}{@{}ccllllllll@{}}
\toprule
\multirow{2}{*}{\textsc{Dataset}} &
  \multirow{2}{*}{\begin{tabular}[c]{@{}c@{}}\textsc{Reconstruction}\\ \textsc{Scheme}\end{tabular}} &
  \multicolumn{2}{c}{2x} &
  \multicolumn{2}{c}{4x} &
  \multicolumn{2}{c}{6x} &
  \multicolumn{2}{c}{8x} \\ \cmidrule(l){3-10} 
                              &              & PSNR           & SSIM           & PSNR           & SSIM           & PSNR           & SSIM  & PSNR          & SSIM  \\ \midrule
\multirow{6}{*}{\textsc{Cardiac Cine}} & Zero-Filling & 31.40          & 90.19          & 25.14          & 81.40          & 25             & 76.97 & 22.99         & 72.71 \\
                              & CS           & 32.85          & 93.58          & 31.57          & 88.57          & 27.98          & 81.48 & 22.97         & 72.76 \\
                              & L+S          & 34.21          & 92.77          & 31.09          & 86.20          & 27.74          & 80.10 & 22.98         & 72.68 \\
                              & CS+M         & 36.72          & 96.23          & 31.53          & 90.26          & \textbf{28.39} & 80.17 & \textbf{24.9} & 72.68 \\
                              & MC+JPDAL     & 36.7           & \textbf{97.85} & \textbf{32.72} & 92.06          & 27.80          & 84.29 & 23.15         & 75.51 \\
                              & DLMCR        & \textbf{38.01} & 97.33          & 32.35          & \textbf{92.26} & 27.65          & 84.76 & 23.12         & 76.03 \\ \hline
\multicolumn{1}{l}{} &
  \begin{tabular}[c]{@{}c@{}} DLMCR w/TL\end{tabular} &
  37 &
  97.05 &
  31.79 &
  90.71 &
  27.70 &
  \textbf{84.85} &
  24.80 &
  \textbf{81.01} \\ \bottomrule
\end{tabular}
}
\caption{Numerical comparison of our technique vs other reconstruction methods. The numerical values are computed as the averages of the similarity metrics over the complete corresponding dataset. w/TL denotes the transfer learning capability of our technique, that is- the results are from training our dictionary with phantom datasets and applied to the real cardiac cine.}
\label{tab:phantom_fb2cine_numerical_results}
\end{table}

%% file: sections/5_conclusion.tex
This work addresses the problem of reconstructing high-quality images from dynamic under-sampled MRI measurements. This is a fundamental problem in MRI due to the long acquisition time required to form an image. During scanning, involuntary and breathing motion often lead to image degradation, compromising clinical interpretation and relevance for diagnosis. 

To tackle this problem, we propose a novel mathematical model to improve the reconstruction quality by estimating motion in dynamic MRI. The underlying idea is to exploit the strong correlation between the two tasks of image reconstruction and motion estimation.The main motivation of considering motion as a second task is because, in a dynamic setting, inherent motion is contained in the scene and therefore it strongly depends on the reconstruction quality. Our research hypothesis is that the reconstruction accuracy can be substantially improved by a good approximation of the motion scene. Our proposed approach combines, in a multi-task and hybrid model, the traditional compressed sensing formulation for the reconstruction of dynamic MRI and motion compensation by learnt optical flow approximation. First, we introduce the classical compressed sensing reconstruction formulation along with motion estimation algorithms depicted in \cite{aviles2018}. Second, we show how the optical flow can be learnt from reference samples following the work introduced in \cite{jia2011}. More precisely, our optical flow model is computed as a sparse linear combination of basis functions from a learnt dictionary. Then, we embed in a single functional the reconstruction process and the optical flow estimation by dictionary learning. 
We present in details our efficient and tractable optimisation framework adopted to solve the non-convex problem based on an alternating splitting scheme. Finally, we extensively demonstrate the potential of our proposed model in the context of dynamic MRI through a set of numerical experiments. In particular, we present comparisons with single task MRI reconstructions as well as with joint MRI reconstruction and motion estimation models, demonstrating that our approach reaches and outperforms state-of-the-art among variational models. 

Our multi-task hybrid model has shown enormous potential for improving the quality of dynamic MRI reconstruction from highly undersampled data by carefully intertwining two imaging tasks that are traditionally performed separately. Furthermore, our motion estimation is based on a learnt approach that relies on a sparse representation in a dictionary of reference samples, with the great advantage that ground truth data is not required. In addition, we explore the potential of our proposed method in the context of transfer learning, that is we learn our dictionary on phantom data and apply it to real dataset. Our results are promising especially in the medical domain, where many times we deal with a severe lack of ground truth data. A perspective for future research is to embed the proposed algorithm in a coarse to fine pyramidal approach. We would also add to the optical flow constraint another term coding the possible brightness variations in the sequence. Finally, we will take our findings and apply them to a more clinical study. Overall, our technique has three major advantages over other techniques. Firstly, our technique statistically improves reconstruction quality over the compared techniques. Secondly, our model requires less computational resources than other multi-task models. Finally, we show the ability of our technique to transfer knowledge from phantom to real data.

%% file: sections/6_appendices.tex
\section{Description of the data structures and the operators}\label{anx:operators}
We assume that the space-time discrete grid consists of the following set of points: 
\begin{align*}
    \left\{(i,j,t)\ :\ i=1,\ldots,N_x,\ j=1,\ldots,N_y,\ t=1,\ldots,N_t\right\}.
\end{align*}
\noindent
Therefore, the reconstructed image sequence $\mathbf{m}$ takes values in $\mathbb{R}^{N_xN_yN_t}$, whereas the components $u_x,\ u_y$ of the optical flow take values in $\mathbb{R}^{N_xN_y(N_t-1)}$. We now discuss the discretisation of the operators. Assume we have a sequence $\mathbf{m}$ of $N_t$ frames of size $N_x\times N_y$ pixels. We use a forward discretisation for the temporal derivative and a central discretisation for the spatial derivative. Neumann boundary conditions are applied: 

\begin{align*}
    \mathbf{m}_t(i,j,t) &= 
        \left\{\begin{array}{ll}
            \mathbf{m}(i,j,t+1) - \mathbf{m}(i,j,t) & \text{if}\ t < N_t,  \\
            0 & \text{else}, 
        \end{array} \right. \\
    \mathbf{m}_x(i,j,t) &= \frac{1}{2}
        \left\{ \begin{array}{ll}
            \mathbf{m}(i+1,j,t) - \mathbf{m}(i-1,j,t) & \text{if}\ i>1\ \text{and}\ i<N_x\ \text{and}\ t < N_t,  \\
            0 & \text{else}, 
        \end{array} \right. \\    
    \mathbf{m}_y(i,j,t) &= \frac{1}{2}
        \left\{ \begin{array}{ll}
            \mathbf{m}(i,j+1,t) - \mathbf{m}(i,j-1,t) & \text{if}\ j>1\ \text{and}\ j<N_y\ \text{and}\ t < N_t,  \\
            0 & \text{else}.
        \end{array} \right. 
\end{align*}

The adjoint operator then yields: 
\begin{align*}
    \mathbf{y}_t(i,j,t) &= -
    \left\{ \begin{array}{ll}
        \mathbf{y}(i,j,t) & \text{if}\ t=1, \\ 
        \mathbf{y}(i,j,t)-\mathbf{y}(i,j,t-1) & \text{if}\ t>1\ \text{and}\ t<N_t, \\ 
        -\mathbf{y}(i,j,t-1) & \text{if}\ t=N_t, \\ 
        \end{array} \right. \\
    \mathbf{y}_x(i,j,t) &= \frac{-1}{2}
    \left\{ \begin{array}{ll}
        \mathbf{y}(i+1,j,y) & \text{if}\ i\leq 2\ \text{and}\ t<N_t \\
        \mathbf{y}(i-1,j,y)-\mathbf{y}(i+1,j,y) & \text{if}\ i>2\ \text{and}\ i< N_x-1\ \text{and}\ t<N_t \\
        -\mathbf{y}(i-1,j,y) & \text{if}\ i\geq N_x-1\ \text{and}\ t<N_t \\
        0 & \text{else},
        \end{array} \right. \\    
    \mathbf{y}_y(i,j,t) &= \frac{-1}{2}
    \left\{ \begin{array}{ll}
        \mathbf{y}(i,j+1,y) & \text{if}\ j\leq 2\ \text{and}\ t<N_t \\
        \mathbf{y}(i,j-1,y)-\mathbf{y}(i,j+1,y) & \text{if}\ j>2\ \text{and}\ j< N_y-1\ \text{and}\ t<N_t \\
        -\mathbf{y}(i,j-1,y) & \text{if}\ j\geq N_y-1\ \text{and}\ t<N_t \\
        0 & \text{else}.
    \end{array} \right. 
\end{align*}

The discrete derivatives for the optical flow vector field are calculated using forward differences and Neumann boundary conditions. The corresponding adjoint operator consists of backward differences with Dirichlet boundary conditions and is applied to the dual variable $\mathbf{y}$. In the following, $u$ denotes either the horizontal component $\mathbf{u}_x$ or the vertical $\mathbf{u}_y$.

\begin{align*}
    u_x(i,j,t) &=
        \left\{ \begin{array}{ll}
            u(i+1,j,t) - u(i,j,t) & \text{if}\ i<N_x \\
            0 & \text{if}\ i=N_x, 
        \end{array} \right. \\    
    u_y(i,j,t) &=
        \left\{ \begin{array}{ll}
            u(i,j+1,t) - u(i,j,t) & \text{if}\ j<N_y \\
            0 & \text{if}\ j=N_y, 
        \end{array} \right. \\    
    \nabla\cdot\mathbf{y}(i,j,t) &= 
        \left\{ \begin{array}{ll}
            \mathbf{y}^1(i,j,t) & \text{if}\ i=1 \\
            \mathbf{y}^1(i,j,t) -\mathbf{y}^1(i-1,j,t) &\text{if}\ i>1\ \text{and}\ i<N_x \\
            -\mathbf{y}^1(i-1,j,t) &\text{if}\ i=N_x
        \end{array} \right.  \\ 
        &\qquad +
        \left\{ \begin{array}{ll}
        \mathbf{y}^2(i,j,t) & \text{if}\ j=1 \\
        \mathbf{y}^2(i,j,t) -\mathbf{y}^2(i,j-1,t) &\text{if}\ j>1\ \text{and}\ j<N_y \\
        -\mathbf{y}^2(i,j-1,t) &\text{if} j=N_y.
        \end{array} \right.
\end{align*}

\section{Error measures}\label{anx:err_measures}
We evaluate the performance of our model in terms of quality for the reconstructed images. 
The main metric thus used is 
the structural similarity index (SSIM) \cite{Wang::2004} which measures the differences in luminance, contrast and structure of the ground truth image $m$ and the reconstruction $m_r$ as follows, 

\begin{align*}
    SSIM := \frac{(2\mu_m\mu_{m_r} + C_1)(2\sigma_{m,m_r} + C_2)}{(\mu_m^2  + \mu_{m_r}^2 + C_1)(\sigma_m^2 + \sigma_{m_r}^2 + C_2)}
\end{align*}

where $\mu_m,\mu_{m_r},\sigma_m,\sigma_{m_r}$ and $\sigma_{m,m_r}$ are local means, standard deviations and cross covariance for ground truth image $m$ and reconstruction $m_r$ respectively. The constants are fixed to $C_1=0.01^2$ and $C_2=0.03^2$. The SSIM takes values between $-1$ and $1$ where $1$ stands for perfect similarity. Moreover we calculate the peak signal-to-noise ratio (PSNR) between the ground truth and the reconstruction given by 

\begin{align}
    PSNR := 10\log_{10}\left(\frac{\max(m^2\textcolor{red}{)}}{\text{mean}((m-m_r)^2)}\right).
\end{align}

\section{Proximal operators}\label{anx:proximal}
In this section we present the computation of the proximal operator of the optical flow minimisation for the Chambolle and Pock algorithm \cite{Chambolle::2011}. We recall the unconstrained problem, 

\begin{align*}
    \min_{\mathbf{u}} \int_0^T \lambda_3\left\|\frac{\partial \mathbf{m}}{\partial t} + \nabla\mathbf{m}\cdot\mathbf{u}\right\|_1 +\lambda_4\|\nabla\mathbf{u}\|_1\  +  \lambda_5 \sum_{\mathbf{p}\in\mathcal{P}} \left\|R_{\mathbf{p}}\mathbf{u}-D\mathbf{a}_{\mathbf{p}}\right\|_F^2 dt .
\end{align*}

To simplify the notation, we introduce the affine operator $\rho(\mathbf{u})=\frac{\partial \mathbf{m}}{\partial t}+\nabla\mathbf{m}\cdot \mathbf{u}$. The challenging part is to compute the resolvent operator of $G(\mathbf{u}) = \lambda_3\|\rho(\mathbf{u})\|_1+\lambda_5 \sum_{\mathbf{p}\in\mathcal{P}}\|R_{\mathbf{p}} \mathbf{u}-$ $D \mathbf{a}_{\mathbf{p}} \|_F^2 $. We have, 

\begin{align}
    \label{eq:anx_prox_of}
    (I+\tau\partial G)^{-1} (\mathbf{u}) &= \arg\min_{\mathbf{v}}\left\{\frac{\|\mathbf{v}-\mathbf{u}\|^2}{2\tau}+ G(\mathbf{v})\right\} \\
\nonumber    &= \arg\min_{\mathbf{v}}\left\{\frac{\|\mathbf{v}-\mathbf{u}\|^2}{2\tau}+\lambda_3\|\rho(\mathbf{v})\|_1 + \lambda_5\sum_{\mathbf{p}\in\mathcal{P}}\|R_{\mathbf{p}}\mathbf{v} - D\mathbf{a}_{\mathbf{p}}\|_F^2\right\}. 
\end{align}

\noindent
The problem is decomposed in three cases: $\rho(\mathbf{u}) > 0$, $\rho(\mathbf{u}) < 0$ and $\rho(\mathbf{u}) = 0$, the previous constraint are intended component-wise. In the first two cases the $l^1$-norm is differentiable so they can be treated in a similar fashion:

\begin{align}
    \rho(\mathbf{v}) >0 &\Longrightarrow \frac{\mathbf{v} - \mathbf{u}}{\tau}+ \lambda_3\nabla\mathbf{m}+2\lambda_5\sum R_{\mathbf{p}}^T (R_{\mathbf{p}}\mathbf{v}-D\mathbf{a}_{\mathbf{p}}) = 0 \nonumber \\ 
\label{eq:anx_prox_of_p_v}    & \Longrightarrow \left(I+2\tau\lambda_5\sum R_{\mathbf{p}}^T R_{\mathbf{p}}\right)\mathbf{v} = \mathbf{u}-\tau\lambda_3\nabla\mathbf{m}+2\tau\lambda_5\sum R_{\mathbf{p}}^TD\mathbf{a}_{\mathbf{p}}.
\end{align}

The equality \cref{eq:anx_prox_of_p_v} gives the solution of \cref{eq:anx_prox_of}. Like in \cref{eq:optflow_prox} we note the operator $\left(I+2\tau\lambda_5\sum R_{\mathbf{p}}^T R_{\mathbf{p}}\right)$ by $\mathcal{A}$. We then plug the solution \cref{eq:anx_prox_of_p_v} in the condition $\rho(\mathbf{v}) >0$ to obtain a necessary condition for the realisation of this solution: 

\begin{align}
    \rho\left(\mathcal{A}^{-1}(\mathbf{u} + 2\tau\lambda_5\sum R_{\mathbf{p}}^TD\mathbf{a}_{\mathbf{p}})\right) - \tau\lambda_3\mathcal{A}^{-1}\|\nabla\mathbf{m}\|^2 > 0.
\end{align}

When $\rho(\mathbf{v})=0$ the $l^1$-norm is no longer differentiable so we have to consider its sub-differential. 

\begin{align}
    \rho(\mathbf{v}) =0 &\Longrightarrow \frac{\mathbf{v} - \mathbf{u}}{\tau}+ \lambda_3 [-\nabla\mathbf{m};\nabla\mathbf{m}] +2\lambda_5\sum R_{\mathbf{p}}^T(R_{\mathbf{p}}\mathbf{v}-D\mathbf{a}_{\mathbf{p}}) = 0 \nonumber \\ 
    &\Longrightarrow \left(I+2\tau\lambda_5\sum R_{\mathbf{p}}^TR_{\mathbf{p}}\right)\mathbf{v} = \mathbf{u} - \tau\lambda_3[-\nabla\mathbf{m};\nabla\mathbf{m}] +2\tau\lambda_5\sum R_{\mathbf{p}}^TD\mathbf{a}_{\mathbf{p}}. \label{eq:anx_prox_of_0}
\end{align}

In the previous \cref{eq:anx_prox_of_0} we note by $[-\nabla\mathbf{m};\nabla\mathbf{m}]$ any element of the sub-differential of the $l^1$-norm, then the condition $\rho(\mathbf{v})=0$ gives us the required element: 

\begin{align}
    \rho\left(\mathcal{A}^{-1}(\mathbf{u} + 2\tau\lambda_5\sum R_{\mathbf{p}}^TD\mathbf{a}_{\mathbf{p}})\right) -\tau\lambda_3[-\nabla\mathbf{m};\nabla\mathbf{m}]\cdot\nabla\mathbf{m}  &=0 \nonumber\\ 
    \Longrightarrow \tau\lambda_3[-\nabla\mathbf{m};\nabla\mathbf{m}] = \rho\left(\mathcal{A}^{-1}(\mathbf{u} + 2\tau\lambda_5\sum R_{\mathbf{p}}^TD\mathbf{a}_{\mathbf{p}})\right) \frac{\nabla\mathbf{m}}{\|\nabla\mathbf{m}\|^2}.
\end{align}

Finally,
\begin{align*}
    \mathcal{A}(I+\tau\partial G)^{-1} (\mathbf{u}) &=\tilde{\mathbf{u}} + \left\{
    \begin{array}{l}
        -\tau\lambda_3\nabla\mathbf{m} \quad \text{if} \quad \rho\left(\mathcal{A}^{-1}\tilde{\mathbf{u}}\right) > \tau\lambda_3\mathcal{A}^{-1}\|\nabla \mathbf{m}\|^2 \\ 
        \tau \lambda_3\nabla\mathbf{m} \quad \text{if} \quad \rho\left(\mathcal{A}^{-1}\tilde{\mathbf{u}}\right) < -\tau\lambda_3\mathcal{A}^{-1}\|\nabla \mathbf{m}\|^2  \\ 
        -\tau\lambda_3\rho\left(\mathcal{A}^{-1}\tilde{\mathbf{u}}\right)\frac{\nabla\mathbf{m}}{\|\nabla\mathbf{m}\|^2} \quad \text{otherwise}
    \end{array}\right.
\end{align*}

where $\tilde{\mathbf{u}} = \mathbf{u} + 2\tau\lambda_5\sum R_{\mathbf{p}}^TD\mathbf{a}_{\mathbf{p}}  $.

\section{Supplementary Experiments} In this section, we extend the comparison results of Table 2 from the main paper. Table~\ref{tab:numerical_results_supplementary} shows a performance comparison of our technique under different initialisations (ZF, CS and CS+M). A closer look at the results shows that as better the initialisation as higher the reconstruction quality. Whilst using initialisation such as CS or CS+M indeed improves the performance, it also implies several drawbacks such as substantial computational load as additional optimisation process need to be solved, and biased reconstructions closer to the initialisation. Moreover, the problem at hand might be seen as being solved twice. For these reasons, we follow common protocol for initialisation and, in the main paper, we use zero filling as our initialisation. 

\medskip

\begin{table}[h]
\centering
\resizebox{1\columnwidth}{!}{
\begin{tabular}{c|cll|cll|cll}
\begin{tabular}[c]{@{}c@{}}\textsc{Acceleration}\\ \textsc{Factor}\end{tabular} & \textsc{Initialisation} & \multicolumn{1}{c}{SSIM} & \multicolumn{1}{c|}{PSNR} & \textsc{Initialisation} & \multicolumn{1}{c}{SSIM} & \multicolumn{1}{c|}{PSNR} & \textsc{Initialisation} & \multicolumn{1}{c}{SSIM} & \multicolumn{1}{c}{PSNR} \\ \hline
2x & \multirow{4}{*}{ZF} & 97.33 & 38.01 & \multirow{4}{*}{CS} & 97.70 & 37.47 & \multirow{4}{*}{CS+M} & 98.34 & 37.75 \\
4x &  & 92.26 & 32.35 &  & 92.47 & 32.37 &  & 92.90 & 33.24 \\
6x &  & 84.76 & 27.65 &  & 84.68 & 28.07 &  & 84.92 & 29.05 \\
8x &  & 76.03 & 23.12 &  & 76.19 & 23.66 &  & 76.61 & 24.40 \\ \hline
\end{tabular}}
\caption{Performance comparison of our technique under different initialisations and different undersampling factors. The results are reported as the average of the corresponding metric on Dataset I.}
\label{tab:numerical_results_supplementary}
\end{table}

%% file: article_color.bbl
\begin{thebibliography}{10}

\bibitem{arif2019accelerated}
{\sc O.~Arif, H.~Afzal, H.~Abbas, M.~F. Amjad, J.~Wan, and R.~Nawaz}, {\em
  {Accelerated Dynamic MRI Using Kernel-Based Low Rank Constraint}}, Journal of
  medical systems, 43 (2019), p.~271.

\bibitem{Asif::2013}
{\sc M.~S. Asif, L.~Hamilton, M.~Brummer, and J.~Romberg}, {\em Motion-adaptive
  spatio-temporal regularization for accelerated dynamic {MRI}}, Magnetic
  Resonance in Medicine, 70 (2013), pp.~800--812.

\bibitem{aviles2018}
{\sc A.~I. Aviles-Rivero, G.~Williams, M.~J. Graves, and C.-B. Schonlieb}, {\em
  Compressed sensing plus motion ({CS+ M}): a new perspective for improving
  undersampled {MR} image reconstruction}, arXiv preprint arXiv:1810.10828,
  (2018).

\bibitem{Aviles-Rivero::2018}
{\sc G.~M. Aviles-Rivero~AI, Williams~G and S.~CB.}, {\em {CS+M}: A
  simultaneous reconstruction and motion estimation approach for improving
  undersampled {MRI} reconstruction}, in In Proceedings of the 26th Annual
  Meeting ISMRM, 2018.

\bibitem{bao2015dictionary}
{\sc C.~Bao, H.~Ji, Y.~Quan, and Z.~Shen}, {\em Dictionary learning for sparse
  coding: {A}lgorithms and convergence analysis}, IEEE transactions on pattern
  analysis and machine intelligence, 38 (2015), pp.~1356--1369.

\bibitem{burger2018}
{\sc M.~Burger, H.~Dirks, and C.~Sch{\"o}nlieb}, {\em A variational model for
  joint motion estimation and image reconstruction}, SIAM Journal on Imaging
  Sciences, 11 (2018), pp.~94--128.

\bibitem{Candes::2011}
{\sc E.~J. Cand{\`e}s, X.~Li, Y.~Ma, and J.~Wright}, {\em Robust principal
  component analysis?}, Journal of the ACM (JACM), 58 (2011), p.~11.

\bibitem{candes2006stable}
{\sc E.~J. Candes, J.~K. Romberg, and T.~Tao}, {\em Stable signal recovery from
  incomplete and inaccurate measurements}, Communications on Pure and Applied
  Mathematics: A Journal Issued by the Courant Institute of Mathematical
  Sciences, 59 (2006), pp.~1207--1223.

\bibitem{caruana1997multitask}
{\sc R.~Caruana}, {\em Multitask learning}, Machine learning, 28 (1997),
  pp.~41--75.

\bibitem{Chambolle::2011}
{\sc A.~Chambolle and T.~Pock}, {\em A first-order primal-dual algorithm for
  convex problems with applications to imaging}, Journal of Mathematical
  Imaging and Vision,  (2011), pp.~120--145.

\bibitem{Corona::2019}
{\sc V.~Corona, A.~Aviles-Rivero, N.~Debroux, M.Grave, C.~L. Guyader, C.-B.
  Sch{\"o}nlieb, and G.~Williams}, {\em Multi-tasking to correct:
  motion-compensated {MRI} via joint reconstruction and registration}, Scale
  Space and Variational Methods in Computer Vision LNCS conference proceedings,
  Springer.,  (2019).

\bibitem{donoho2006compressed}
{\sc D.~L. Donoho}, {\em Compressed sensing}, IEEE Transactions on information
  theory, 52 (2006), pp.~1289--1306.

\bibitem{dosovitskiy2015flownet}
{\sc A.~Dosovitskiy, P.~Fischer, E.~Ilg, P.~Hausser, C.~Hazirbas, V.~Golkov,
  P.~Van Der~Smagt, D.~Cremers, and T.~Brox}, {\em Flownet: Learning optical
  flow with convolutional networks}, in IEEE international conference on
  computer vision (CVPR), 2015, pp.~2758--2766.

\bibitem{elad2006image}
{\sc M.~Elad and M.~Aharon}, {\em Image denoising via sparse and redundant
  representations over learned dictionaries}, IEEE Transactions on Image
  Processing, 15 (2006), pp.~3736--3745.

\bibitem{Feng::2013}
{\sc L.~Feng, M.~B. Srichai, R.~P. Lim, A.~Harrison, W.~King, G.~Adluru, E.~V.
  Dibella, D.~K. Sodickson, R.~Otazo, and D.~Kim}, {\em Highly accelerated
  real-time cardiac cine {MRI} using k--t {SPARSE-SENSE}}, Magnetic Resonance
  in Medicine,  (2013), pp.~64--74.

\bibitem{filipovic2011motion}
{\sc M.~Filipovic, P.-A. Vuissoz, A.~Codreanu, M.~Claudon, and J.~Felblinger},
  {\em Motion compensated generalized reconstruction for free-breathing dynamic
  contrast-enhanced {MRI}}, Magnetic Resonance in Medicine, 65 (2011),
  pp.~812--822.

\bibitem{fortun2015}
{\sc D.~Fortun, P.~Bouthemy, and C.~Kervrann}, {\em Optical flow modeling and
  computation: A survey}, Computer Vision and Image Understanding, 134 (2015),
  pp.~1--21.

\bibitem{gamper2008compressed}
{\sc U.~Gamper, P.~Boesiger, and S.~Kozerke}, {\em Compressed sensing in
  dynamic {MRI}}, Magnetic Resonance in Medicine: An Official Journal of the
  International Society for Magnetic Resonance in Medicine, 59 (2008),
  pp.~365--373.

\bibitem{Gao::2012}
{\sc H.~Gao, S.~Rapacchi, D.~Wang, J.~Moriarty, C.~Meehan, J.~Sayre, G.~Laub,
  P.~Finn, and P.~Hu}, {\em Compressed sensing using prior rank, intensity and
  sparsity model ({PRISM}): applications in cardiac cine {MRI}}, in Proceedings
  of the Annual Meeting of ISMRM, 2012.

\bibitem{gdaniec2014robust}
{\sc N.~Gdaniec, H.~Eggers, P.~B{\"o}rnert, M.~Doneva, and A.~Mertins}, {\em
  Robust abdominal imaging with incomplete breath-holds}, Magnetic Resonance in
  Medicine, 71 (2014), pp.~1733--1742.

\bibitem{GEORGE2006924}
{\sc R.~George, T.~D. Chung, S.~S. Vedam, V.~Ramakrishnan, R.~Mohan, E.~Weiss,
  and P.~J. Keall}, {\em Audio-visual biofeedback for respiratory-gated
  radiotherapy: Impact of audio instruction and audio-visual biofeedback on
  respiratory-gated radiotherapy}, International Journal of Radiation
  Oncology*Biology*Physics, 65 (2006), pp.~924 -- 933.

\bibitem{gibson2016sparsity}
{\sc J.~Gibson and O.~Marques}, {\em Sparsity in optical flow and
  trajectories}, Signal, Image and Video Processing, 10 (2016), pp.~487--494.

\bibitem{haskell2018targeted}
{\sc M.~W. Haskell, S.~F. Cauley, and L.~L. Wald}, {\em {TArgeted Motion
  Estimation and Reduction (TAMER): data consistency based motion mitigation
  for MRI using a reduced model joint optimization}}, IEEE transactions on
  medical imaging, 37 (2018), pp.~1253--1265.

\bibitem{HOISAK2006339}
{\sc J.~D. Hoisak, K.~E. Sixel, R.~Tirona, P.~C. Cheung, and J.-P. Pignol},
  {\em Prediction of lung tumour position based on spirometry and on abdominal
  displacement: Accuracy and reproducibility}, Radiotherapy and Oncology, 78
  (2006), pp.~339 -- 346.

\bibitem{horn1981}
{\sc B.~K. Horn and B.~G. Schunck}, {\em Determining optical flow}, in
  Techniques and Applications of Image Understanding, vol.~281, International
  Society for Optics and Photonics, 1981, pp.~319--331.

\bibitem{ilg2017flownet}
{\sc E.~Ilg, N.~Mayer, T.~Saikia, M.~Keuper, A.~Dosovitskiy, and T.~Brox}, {\em
  Flownet 2.0: Evolution of optical flow estimation with deep networks}, in
  IEEE Conference on Computer Vision and Pattern Recognition (CVPR), 2017,
  pp.~2462--2470.

\bibitem{jia2011}
{\sc K.~Jia, X.~Wang, and X.~Tang}, {\em Optical flow estimation using learned
  sparse model}, in 2011 International Conference on Computer Vision, IEEE,
  2011, pp.~2391--2398.

\bibitem{JIANG2006141}
{\sc S.~B. Jiang}, {\em Technical aspects of image-guided respiration-gated
  radiation therapy}, Medical Dosimetry, 31 (2006), pp.~141 -- 151.

\bibitem{jung2009k}
{\sc H.~Jung, K.~Sung, K.~S. Nayak, E.~Y. Kim, and J.~C. Ye}, {\em k-t
  {FOCUSS}: a general compressed sensing framework for high resolution dynamic
  {MRI}}, Magnetic Resonance in Medicine: An Official Journal of the
  International Society for Magnetic Resonance in Medicine, 61 (2009),
  pp.~103--116.

\bibitem{Kaji2001}
{\sc S.~Kaji, P.~C. Yang, A.~B. Kerr, W.~W. Tang, C.~H. Meyer, A.~Macovski,
  J.~M. Pauly, D.~G. Nishimura, and B.~S. Hu}, {\em Rapid evaluation of left
  ventricular volume and mass without breath-holding using real-time
  interactive cardiac magnetic resonance imaging system}, Journal of the
  American College of Cardiology, 38 (2001), pp.~527--533.

\bibitem{li2019rainflow}
{\sc R.~Li, R.~T. Tan, L.-F. Cheong, A.~I. Aviles-Rivero, Q.~Fan, and C.-B.
  Sch{\"o}nlieb}, {\em {RainFlow}: Optical flow under rain streaks and rain
  veiling effect}, in IEEE International Conference on Computer Vision (ICCV),
  2019, pp.~7304--7313.

\bibitem{Liang::2007}
{\sc Z.-P. Liang}, {\em Spatiotemporal imaging with partially separable
  functions}, in IEEE International Symposium on Biomedical Imaging (ISBI),
  2007, pp.~988--991.

\bibitem{Lingala::2015}
{\sc S.~G. Lingala, E.~DiBella, and M.~Jacob}, {\em Deformation corrected
  compressed sensing ({DC-CS}): a novel framework for accelerated dynamic
  {MRI}}, IEEE transactions on medical imaging, 34 (2014), pp.~72--85.

\bibitem{Lingala::2011}
{\sc S.~G. Lingala, Y.~Hu, E.~DiBella, and M.~Jacob}, {\em Accelerated dynamic
  {MRI} exploiting sparsity and low-rank structure: kt {SLR}}, IEEE
  Transactions on Medical Imaging,  (2011), pp.~1042--1054.

\bibitem{lingala::2013}
{\sc S.~G. Lingala and M.~Jacob}, {\em Blind compressive sensing dynamic mri},
  IEEE transactions on medical imaging, 32 (2013), pp.~1132--1145.

\bibitem{liu2019rethinking}
{\sc J.~Liu, A.~I. Aviles-Rivero, H.~Ji, and C.-B. Sch{\"o}nlieb}, {\em
  Rethinking medical image reconstruction via shape prior, going deeper and
  faster: Deep joint indirect registration and reconstruction}, arXiv preprint
  arXiv:1912.07648,  (2019).

\bibitem{lustig2007}
{\sc M.~Lustig, D.~Donoho, and J.~M. Pauly}, {\em Sparse {MRI}: The application
  of compressed sensing for rapid {MR} imaging}, Magnetic Resonance in
  Medicine, 58 (2007), pp.~1182--1195.

\bibitem{LustigSPIRiT::2010}
{\sc M.~Lustig and J.~M. Pauly}, {\em {SPIRiT}: iterative self-consistent
  parallel imaging reconstruction from arbitrary k-space}, Magnetic Resonance
  in Medicine,  (2010), pp.~457--471.

\bibitem{LustigktSPARSE::2006}
{\sc M.~Lustig, J.~M. Santos, D.~L. Donoho, and J.~M. Pauly}, {\em k-t
  {SPARSE}: High frame rate dynamic {MRI} exploiting spatio-temporal sparsity},
  in Proceedings of the Annual Meeting of ISMRM, vol.~2420, 2006.

\bibitem{mairal2009non}
{\sc J.~Mairal, F.~Bach, J.~Ponce, G.~Sapiro, and A.~Zisserman}, {\em Non-local
  sparse models for image restoration}, in 2009 IEEE 12th International
  Conference on Computer Vision, IEEE, 2009, pp.~2272--2279.

\bibitem{Majumdar::2015}
{\sc A.~Majumdar}, {\em Compressed sensing for magnetic resonance image
  reconstruction}, Cambridge University Press, 2015.

\bibitem{Majumdar::2012}
{\sc A.~Majumdar and R.~K. Ward}, {\em Exploiting rank deficiency and transform
  domain sparsity for {MR} image reconstruction}, Magnetic Resonance Imaging,
  (2012), pp.~9--18.

\bibitem{Miao::2016}
{\sc X.~Miao, S.~G. Lingala, Y.~Guo, T.~Jao, M.~Usman, C.~Prieto, and K.~S.
  Nayak}, {\em Accelerated cardiac cine {MRI} using locally low rank and finite
  difference constraints}, Magnetic Resonance Imaging,  (2016).

\bibitem{Odille::2016}
{\sc F.~Odille, A.~Menini, J.-M. Escany{\'e}, P.-A. Vuissoz, P.-Y. Marie,
  M.~Beaumont, and J.~Felblinger}, {\em Joint reconstruction of multiple images
  and motion in {MRI}: Application to free-breathing myocardial ${\rm{t}}_{2}$
  quantification}, IEEE transactions on medical imaging, 35 (2015),
  pp.~197--207.

\bibitem{Odille::2010}
{\sc F.~Odille, S.~Uribe, P.~G. Batchelor, C.~Prieto, T.~Schaeffter, and
  D.~Atkinson}, {\em Model-based reconstruction for cardiac cine {MRI} without
  {ECG} or breath holding}, Magnetic Resonance in Medicine, 63 (2010),
  pp.~1247--1257.

\bibitem{Odille::2008}
{\sc F.~Odille, P.-A. Vuissoz, P.-Y. Marie, and J.~Felblinger}, {\em
  Generalized reconstruction by inversion of coupled systems ({GRICS}) applied
  to free-breathing {MRI}}, Magnetic Resonance in Medicine: An Official Journal
  of the International Society for Magnetic Resonance in Medicine, 60 (2008),
  pp.~146--157.

\bibitem{olshausen1996emergence}
{\sc B.~A. Olshausen and D.~J. Field}, {\em Emergence of simple-cell receptive
  field properties by learning a sparse code for natural images}, Nature, 381
  (1996), pp.~607--609.

\bibitem{Otazo::2015}
{\sc R.~Otazo, E.~Cand{\`e}s, and D.~K. Sodickson}, {\em Low-rank plus sparse
  matrix decomposition for accelerated dynamic {MRI} with separation of
  background and dynamic components}, Magnetic Resonance in Medicine,  (2015),
  pp.~1125--1136.

\bibitem{Otazo::2013}
{\sc R.~Otazo, C.~Emmanuel, and D.~K. Sodickson}, {\em Low-rank \& sparse
  matrix decomposition for accelerated {DCE-MRI} with background \& contrast
  separation}, in Proceedings of the ISMRM Workshop on Data Sampling and Image
  Reconstruction, 2013.

\bibitem{Otazo::2010}
{\sc R.~Otazo, D.~Kim, L.~Axel, and D.~K. Sodickson}, {\em Combination of
  compressed sensing and parallel imaging for highly accelerated first-pass
  cardiac perfusion {MRI}}, Magnetic Resonance in Medicine,  (2010),
  pp.~767--776.

\bibitem{ouzir2018robust}
{\sc N.~Ouzir, A.~Basarab, O.~Lairez, and J.-Y. Tourneret}, {\em Robust optical
  flow estimation in cardiac ultrasound images using a sparse representation},
  IEEE Transactions on Medical Imaging,  (2018), pp.~741--752.

\bibitem{Pedersen::2009}
{\sc H.~Pedersen, S.~Kozerke, S.~Ringgaard, K.~Nehrke, and W.~Y. Kim}, {\em k-t
  {PCA}: temporally constrained k-t {BLAST} reconstruction using principal
  component analysis}, Magnetic Resonance in Medicine,  (2009), pp.~706--716.

\bibitem{perez2013}
{\sc J.~S. P{\'e}rez, E.~Meinhardt-Llopis, and G.~Facciolo}, {\em {TV-L1}
  optical flow estimation}, Image Processing On Line, 2013 (2013),
  pp.~137--150.

\bibitem{plathow2006assessment}
{\sc C.~Plathow, S.~Ley, J.~Zaporozhan, M.~Sch{\"o}binger, E.~Gruenig,
  M.~Puderbach, M.~Eichinger, H.-P. Meinzer, I.~Zuna, and H.-U. Kauczor}, {\em
  Assessment of reproducibility and stability of different breath-hold
  maneuvres by dynamic {MRI}: comparison between healthy adults and patients
  with pulmonary hypertension}, European radiology, 16 (2006), pp.~173--179.

\bibitem{Plein2001}
{\sc S.~Plein, W.~H. Smith, J.~P. Ridgway, A.~Kassner, D.~J. Beacock, T.~N.
  Bloomer, and M.~U. Sivananthan}, {\em Qualitative and quantitative analysis
  of regional left ventricular wall dynamics using real-time magnetic resonance
  imaging: comparison with conventional breath-hold gradient echo acquisition
  in volunteers and patients}, Journal of Magnetic Resonance Imaging: An
  Official Journal of the International Society for Magnetic Resonance in
  Medicine, 14 (2001), pp.~23--30.

\bibitem{Prieto::2007}
{\sc C.~Prieto, P.~G. Batchelor, D.~Hill, J.~V. Hajnal, M.~Guarini, and
  P.~Irarrazaval}, {\em Reconstruction of undersampled dynamic images by
  modeling the motion of object elements}, Magnetic Resonance in Medicine, 57
  (2007), pp.~939--949.

\bibitem{raghu2019transfusion}
{\sc M.~Raghu, C.~Zhang, J.~Kleinberg, and S.~Bengio}, {\em Transfusion:
  Understanding transfer learning for medical imaging}, in Advances in neural
  information processing systems, 2019, pp.~3347--3357.

\bibitem{Rank::2017}
{\sc C.~M. Rank, T.~Heu{\ss}er, M.~T. Buzan, A.~Wetscherek, M.~T. Freitag,
  J.~Dinkel, and M.~Kachelrie{\ss}}, {\em {4D} respiratory motion-compensated
  image reconstruction of free-breathing radial {MR} data with very high
  undersampling}, Magnetic Resonance in Medicine, 77 (2017), pp.~1170--1183.

\bibitem{royuela2016nonrigid}
{\sc J.~Royuela-del Val, L.~Cordero-Grande, F.~Simmross-Wattenberg,
  M.~Mart{\'\i}n-Fern{\'a}ndez, and C.~Alberola-L{\'o}pez}, {\em Nonrigid
  groupwise registration for motion estimation and compensation in compressed
  sensing reconstruction of breath-hold cardiac cine {MRI}}, Magnetic Resonance
  in Medicine, 75 (2016), pp.~1525--1536.

\bibitem{Royuela::2017}
{\sc J.~Royuela-del Val, L.~Cordero-Grande, F.~Simmross-Wattenberg,
  M.~Martín-Fernández, and C.~Alberola-L{\'o}pez}, {\em Jacobian weighted
  temporal total variation for motion compensated compressed sensing
  reconstruction of dynamic {MRI}}, Magnetic Resonance in Medicine, 77 (2017),
  pp.~1208--1215.

\bibitem{sachs1995diminishing}
{\sc T.~S. Sachs, C.~H. Meyer, P.~Irarrazabal, B.~S. Hu, D.~G. Nishimura, and
  A.~Macovski}, {\em The diminishing variance algorithm for real-time reduction
  of motion artifacts in {MRI}}, Magn Reson Med, 34 (1995), pp.~412--422.

\bibitem{Saucedo::2017}
{\sc A.~Saucedo, S.~Lefkimmiatis, N.~Rangwala, and K.~Sung}, {\em Improved
  computational efficiency of locally low rank {MRI} reconstruction using
  iterative random patch adjustments}, IEEE transactions on medical imaging,
  (2017), pp.~1209--1220.

\bibitem{Schloegl::2017}
{\sc M.~Schloegl, M.~Holler, A.~Schwarzl, K.~Bredies, and R.~Stollberger}, {\em
  Infimal convolution of total generalized variation functionals for dynamic
  {MRI}}, Magnetic Resonance in Medicine,  (2017), pp.~142--155.

\bibitem{Setser2000}
{\sc R.~M. Setser, S.~E. Fischer, and C.~H. Lorenz}, {\em Quantification of
  left ventricular function with magnetic resonance images acquired in real
  time}, J Magn Reson Imaging, 12 (2000), pp.~430--438.

\bibitem{shen2010sparsity}
{\sc X.~Shen and Y.~Wu}, {\em Sparsity model for robust optical flow estimation
  at motion discontinuities}, in 2010 IEEE Computer Society Conference on
  Computer Vision and Pattern Recognition, IEEE, 2010, pp.~2456--2463.

\bibitem{sotiras2013}
{\sc A.~{Sotiras}, C.~{Davatzikos}, and N.~{Paragios}}, {\em Deformable medical
  image registration: A survey}, IEEE Transactions on Medical Imaging, 32
  (2013), pp.~1153--1190.

\bibitem{sun2008learning}
{\sc D.~Sun, S.~Roth, J.~Lewis, and M.~J. Black}, {\em Learning optical flow},
  in European Conference on Computer Vision, Springer, 2008, pp.~83--97.

\bibitem{sun2018pwc}
{\sc D.~Sun, X.~Yang, M.-Y. Liu, and J.~Kautz}, {\em {PWC-Net: CNNs for optical
  flow using pyramid, warping, and cost volume}}, in IEEE Conference on
  Computer Vision and Pattern Recognition (CVPR), 2018, pp.~8934--8943.

\bibitem{timofte2015sparse}
{\sc R.~Timofte and L.~Van~Gool}, {\em Sparse flow: Sparse matching for small
  to large displacement optical flow}, in IEEE Winter Conference on
  Applications of Computer Vision (WACV), 2015, pp.~1100--1106.

\bibitem{Tremoulheac::2014}
{\sc B.~Tr{\'e}moulh{\'e}ac, N.~Dikaios, D.~Atkinson, and S.~R. Arridge}, {\em
  {Dynamic {MR} Image Reconstruction--Separation From Undersampled k,t--Space
  via Low-Rank Plus Sparse Prior}}, IEEE Transactions on Medical Imaging,
  (2014), pp.~1689--1701.

\bibitem{Trzasko::2011}
{\sc J.~Trzasko, A.~Manduca, and E.~Borisch}, {\em Local versus global low-rank
  promotion in dynamic {MRI} series reconstruction}, in Proc. Int. Symp. Magn.
  Reson. Med, 2011, p.~4371.

\bibitem{tseng2001convergence}
{\sc P.~Tseng}, {\em Convergence of a block coordinate descent method for
  nondifferentiable minimization}, Journal of optimization theory and
  applications, 109 (2001), pp.~475--494.

\bibitem{Velikina::2015}
{\sc J.~V. Velikina and A.~A. Samsonov}, {\em Reconstruction of dynamic image
  series from undersampled {MRI} data using data-driven model consistency
  condition ({MOCCO})}, Magnetic Resonance in Medicine,  (2015),
  pp.~1279--1290.

\bibitem{von20074d}
{\sc M.~von Siebenthal, G.~Szekely, U.~Gamper, P.~Boesiger, A.~Lomax, and
  P.~Cattin}, {\em {4D MR} imaging of respiratory organ motion and its
  variability}, Physics in Medicine \& Biology, 52 (2007), p.~1547.

\bibitem{Wang::2004}
{\sc Z.~Wang, A.~C. Bovik, H.~R. Sheikh, and E.~P. Simoncelli}, {\em Image
  quality assessment: from error visibility to structural similarity}, IEEE
  Transactions on Image Processing (TIP),  (2004), pp.~600--612.

\bibitem{weller2019motion}
{\sc D.~S. Weller, L.~Wang, J.~P. Mugler~III, and C.~H. Meyer}, {\em
  Motion-compensated reconstruction of magnetic resonance images from
  undersampled data}, Magnetic resonance imaging, 55 (2019), pp.~36--45.

\bibitem{werlberger2010motion}
{\sc M.~Werlberger, T.~Pock, and H.~Bischof}, {\em Motion estimation with
  non-local total variation regularization}, in 2010 IEEE Computer Society
  Conference on Computer Vision and Pattern Recognition, IEEE, 2010,
  pp.~2464--2471.

\bibitem{Wissmann::2014}
{\sc L.~Wissmann, C.~Santelli, W.~P. Segars, and S.~Kozerke}, {\em {MRXCAT}:
  Realistic numerical phantoms for cardiovascular magnetic resonance}, Journal
  of Cardiovascular Magnetic Resonance,  (2014), p.~63.

\bibitem{wulff2015efficient}
{\sc J.~Wulff and M.~J. Black}, {\em Efficient sparse-to-dense optical flow
  estimation using a learned basis and layers}, in Proceedings of the IEEE
  Conference on Computer Vision and Pattern Recognition, 2015, pp.~120--130.

\bibitem{Zaitsev::2015}
{\sc M.~Zaitsev, J.~Maclaren, and M.~Herbst}, {\em Motion artifacts in {MRI}: a
  complex problem with many partial solutions}, Journal of Magnetic Resonance
  Imaging,  (2015), pp.~887--901.

\bibitem{Zhang::2018}
{\sc D.~Zhang, J.~Tao, Z.~Ye, B.~Qiu, and J.~Xu}, {\em Deformation corrected
  blind compressed sensing ({DC-BCS}) a novel framework for dynamic {MRI}
  reconstruction}, in Proceedings of the 2018 2nd International Conference on
  Algorithms, Computing and Systems, 2018, pp.~160--164.

\bibitem{Zhang::2015}
{\sc T.~Zhang, J.~M. Pauly, and I.~R. Levesque}, {\em Accelerating parameter
  mapping with a locally low rank constraint}, Magnetic Resonance in Medicine,
  (2015), pp.~655--661.

\bibitem{Zhao::2012}
{\sc B.~Zhao, J.~P. Haldar, A.~G. Christodoulou, and Z.-P. Liang}, {\em Image
  reconstruction from highly undersampled (k, t)-space data with joint partial
  separability and sparsity constraints}, IEEE Transactions on Medical Imaging,
   (2012), pp.~1809--1820.

\bibitem{Zhao::2019}
{\sc N.~{Zhao}, D.~{O'Connor}, A.~{Basarab}, D.~{Ruan}, and K.~{Sheng}}, {\em
  Motion compensated dynamic {MRI} reconstruction with local affine optical
  flow estimation}, IEEE Transactions on Biomedical Engineering,  (2019),
  pp.~1--1.

\end{thebibliography}
